\documentclass[12pt]{article}
\usepackage{amssymb}
\begin{document}

%Macro for the adjoint:
\def\ad{\mbox{ad}\,}
%Macro for the black-board bold
\def\b#1{{\mathbb #1}}
%Macro for the calligraphic letters
\def\c#1{{\cal #1}}
%Macro to squeeze equation arrays
\def\squeeze{\hspace{-17pt}}
%Macro for the `bra':
\def\bra#1{\langle #1 \vert}
%Macro for the `ket':
\def\ket#1{\vert #1 \rangle}
%Macro for the vacuum expectation value:
\def\vev#1{\langle #1 \rangle}

%Macro to enumerate formulae according to sections and subsections
\newcommand{\sect}[1]{\setcounter{equation}{0}\section{#1}}
\renewcommand{\theequation}{\thesection.\arabic{equation}}
\newcommand{\subsect}[1]{\setcounter{equation}{0}\subsection{#1}}
\renewcommand{\theequation}{\thesection.\arabic{equation}}
%
%%%%%%%%%%%%%%%%%%%%%%%%%%%%%%%
%
\def\1{{\bf 1}}
\newcommand{\be}{\begin{equation}}
\newcommand{\ee}{\end{equation}}
\newcommand{\ba}{\begin{array}}
\newcommand{\ea}{\end{array}}
\newcommand{\bea}{\begin{eqnarray}}
\newcommand{\eea}{\end{eqnarray}}
\newcommand{\beas}{\begin{eqnarray*}}
\newcommand{\eeas}{\end{eqnarray*}}
%
%%%%%%%%%%%%%%%%%%%%%%%%%%%%%%%
%
\newtheorem{prop}{Proposition}
\newtheorem{lemma}{Lemma}
\newtheorem{theorem}{Theorem}
\newtheorem{corollary}{Corollary}
%
%%%%%%%%%%%%%%%%%%%%%%%%%%%%%%%
%

\title{The Geometry of the Quantum Euclidean Space}

\author{Gaetano Fiore$^{1,2}$, \  John Madore$^{3,4}$
        \and
        $\strut^1$Dip. di Matematica e Applicazioni, Fac.  di Ingegneria\\ 
        Universit\`a di Napoli, V. Claudio 21, 80125 Napoli
        \and
        $\strut^2$I.N.F.N., Sezione di Napoli,\\
        Mostra d'Oltremare, Pad. 19, 80125 Napoli
        \and
        $\strut^3$Max-Planck-Institut f\"ur Physik 
        (Werner-Heisenberg-Institut)\\
        F\"ohringer Ring 6, D-80805 M\"unchen
        \and
        $\strut^4$Laboratoire de Physique Th\'eorique et Hautes Energies\\
        Universit\'e de Paris-Sud, B\^atiment 211, F-91405 Orsay
        }

\date{}

\maketitle

\abstract{A detailed study is made of the noncommutative geometry of
$\b{R}^3_q$, the quantum space covariant under the quantum group
$SO_q(3)$.  For each of its two $SO_q(3)$-covariant differential calculi 
we find its metric, the corresponding frame and two
torsion-free covariant derivatives that are metric compatible up to a
conformal factor and which yield both a vanishing linear curvature.  A
discussion is given of various ways of imposing reality conditions. The
delicate issue of the commutative limit is discussed 
at the formal algebraic level.
Two rather different ways of taking the limit are suggested, yielding
respectively $S^2\times \b{R}$ and $\b{R}^3$ as the limit Riemannian
manifold.}

\vfill
\noindent
Preprint LMU-TPW  97-23

\noindent
Preprint 98-63, Dip. Matematica e Applicazioni, Universit\`a di Napoli
\medskip
\eject

\parskip 4pt plus2pt minus2pt

\sect{Introduction and motivation}

Already in 1947 Snyder, in an attempt to remove the divergences from
electrodynamics, suggested~\cite{Sny47} the possibility that the
micro-structure of space-time at the Planck level might be better
described using a noncommutative geometry.  A few years later, in 1954,
Pauli suggested that the gravitational field might be considered as a
universal regularizer for all quantum-field divergences.  An obvious way
to reconciliate these two points of view is to try to define a
gravitational field within the context of noncommutative geometry and
see if any connection between the commutation relations and gravity can
be found. It is our conjecture that in some yet-to-be-found sense the
gravitational field is a classical manifestation of the noncommutative
micro-structure of space-time. Following the usage in ordinary geometry
we shall identify the gravitational field as a metric-compatible,
torsion-free linear connection. We shall use the expression
`connection' and `covariant derivative' synonymously.  

In Section~2 we give a brief overview of a particular version of
noncommutative geometry which can be considered as a noncommutative
extension of the moving-frame formalism if E. Cartan.  More details can
be found elsewhere~\cite{Mad95}. There is not as yet a completely
satisfactory definition of either a linear connection or a metric within
the context of noncommutative geometry but there are definitions which
seem to work in certain cases.  In the present article we add to the
list of examples which seem to lend weight to one particular
definition~\cite{DubMadMasMou96}. We refer to a recent review
article~\cite{Mad97} for a list of some other examples and references to
alternative definitions. More details of one alternative version can be
found in the book by Landi~\cite{Lan97}; for a general introduction to
more mathematical aspects of the subject we refer to the book by
Connes~\cite{Con94}.

The example we consider here is the 3-dim ``quantum Euclidean space''
~\cite{frt}, which has been previously studied by one of the authors.
This is the quantum space covariant under the quantum group
$SO_q(3)$. In Section~3
we give a review of this space. In Section~4 we recall the construction
of the $SO_q(3)$-covariant differential calculi~\cite{cawa} over them,
more specifically the two $SO_q(3)$-covariant calculi determined in the
$\hat R$-matrix formalism \cite{WesZum90} 
and having the de Rham calculus as the commutative limit. 
More details are to be found
in the thesis of one of the authors~\cite{fiothesis}. We shall find that
in order to introduce a `frame' as defined in Section~2 a new generator
for dilatation is needed, as in the construction of inhomogeneous
quantum groups from homogeneous ones, as well as the inverse of some
generator of the quantum Euclidean space.  In Section~5 we define a
metric and two torsion-free linear connection on the quantum Euclidean
space, yielding vanishing linear curvatures.  In Section~6 we consider
the commutative limit $q\to 1$ in order to determine the Riemannian
manifold which remains as a `shadow' of the noncommutative structure.
We suggest two rather different prescriptions for taking the limit: the
limit manifold is $S^2\times \b{R}$ according to the first (and simpler)
prescription, $\b{R}^3$ according to the second (and more sophisticated)
one. The initial generators $x^i$ of the quantum space in the former
tend to their natural cartesian limits, in the latter are to be
`renormalized' before taking the limit.  For $q\neq 1$ $x^i$ appear as a
non-commutative analog of general (i.e. non-cartesian) coordinates. In
Section~7 we discuss the various ways one can construct a real
differential calculus from a combination of the two canonical ones.  In
the last section we present our conclusions and we compare our results
with alternative definitions of linear connections. We compare also our
results with those found in the case of the Manin plane~\cite{DimMad96}
as well as similar results found~\cite{CerHinMadWes98} in the case of
the quantum Euclidean plane of dimension one.

\sect{Linear connections in the frame formalism}               \label{preli}

The starting point is a noncommutative algebra $\c{A}$ which has as
commutative limit the algebra of functions on some manifold and over
$\c{A}$ a differential calculus~\cite{Con94} $\Omega^*(\c{A})$ which
has as corresponding limit the ordinary de~Rham differential
calculus. We recall that a differential calculus is completely
determined by the left and right module structure of the
$\c{A}$-module of 1-forms $\Omega^1(\c{A})$.  We shall restrict our
attention to the case where this module is free of rank $n$ as a left
or right module and possesses a special basis $\theta^a$, 
$1\leq a \leq n$, which commutes with the elements $f$ of the algebra:
\be
[f, \theta^a] = 0.                                             \label{fund} 
\ee
In particular the limit manifold must be parallelizable.  We shall refer
to the $\theta^a$ as a frame or Stehbein.  The integer $n$ plays the
role of `dimension'; it can be greater than the dimension of the limit
manifold but in this case the frame will have a singular limit. We
suppose further~\cite{DimMad96} that the basis is dual to a set of inner
derivations $e_a = \ad \lambda_a$. This means that the differential is
given by the expression
\be
df = e_a f \theta^a = [\lambda_a, f] \theta^a.                 \label{defdiff}
\ee
One can rewrite this equation as 
\be
df = -[\theta,f],                                              \label{extra}
\ee
if one introduces~\cite{Con94} the `Dirac operator'
\be
\theta = - \lambda_a \theta^a.                                 \label{dirac}
\ee

There is a bimodule map $\pi$ of the space 
$\Omega^1(\c{A}) \otimes_\c{A} \Omega^1(\c{A})$ onto the space
$\Omega^2(\c{A})$ of 2-forms and we can write
\be
\theta^a \theta^b = P^{ab}{}_{cd} \theta^c \otimes \theta^d     \label{proj}
\ee
where, because of (\ref{fund}), the $P^{ab}{}_{cd}$ belong to the center
$\c{Z}(\c{A})$ of $\c{A}$. We shall suppose that the center is
trivial, $\c{Z}(\c{A}) = \b{C}$, and therefore the
components $P^{ab}{}_{cd}$ are real numbers. Define the Maurer-Cartan 
elements $C^a{}_{bc} \in \c{A}$ by the equation
\be
d\theta^a = - {1\over 2} C^a{}_{bc} \theta^b \theta^c.       \label{rot-coef}
\ee
Because of~(\ref{proj}) we can suppose that 
$C^a{}_{bc} P^{bc}{}_{de} = C^a{}_{de}$. It follows from the equation
$d(\theta^a f - f \theta^a) = 0$ that there exist elements
$F^a{}_{bc}$ of the center such that
\be
C^a{}_{bc} = F^a{}_{bc} - 2 \lambda_e P^{(ae)}{}_{bc}        \label{consis1}
\ee
where $(ab)$ indicates symmetrization of the indices $a$ and $b$.
If on the other hand we define $K_{ab}$ by the equation
\be
d\theta + \theta^2 = {1\over 2} K_{ab} \theta^a \theta^b,      \label{consis2}
\ee
then if follows from (2.3) and the identity $d^2 = 0$ that the $K_{ab}$
must belong to the center. Finally it can be shown~\cite{DimMad96,
MadMou98} that in order that~(\ref{consis1}) and (\ref{consis2}) be
consistent with one another the original $\lambda_a$ must satisfy the
condition
\be
2 \lambda_c \lambda_d P^{cd}{}_{ab} - 
\lambda_c F^c{}_{ab} - K_{ab} = 0.                             \label{manca}
\ee
This gives to the set of $\lambda_a$ the structure of a twisted Lie
algebra with a central extension.

We propose as definition of a linear connection a map~\cite{Kos60}
\be
\Omega^1(\c{A}) \buildrel D \over \longrightarrow 
\Omega^1(\c{A}) \otimes_\c{A} \Omega^1(\c{A})                   \label{2.2.4} 
\ee 
which satisfies both a left Leibniz rule
\be
D (f \xi) =  df \otimes \xi + f D\xi                            \label{2.2.2}
\ee
and a right Leibniz rule~\cite{DubMadMasMou96}
\be
D(\xi f) = \sigma (\xi \otimes df) + (D\xi) f               %   \label{2.2.3}
\ee
for arbitrary $f \in \c{A}$ and $\xi \in \Omega^1(\c{A})$.  
We have here introduced a generalized flip
\be
\Omega^1(\c{A}) \otimes_\c{A} \Omega^1(\c{A})
\buildrel \sigma \over \longrightarrow
\Omega^1(\c{A}) \otimes_\c{A} \Omega^1(\c{A})                  \label{2.2.5}
\ee
in order to define a right Leibniz rule which is consistent with the
left one. It is necessarily bilinear.  A linear connection is therefore
a couple $(D, \sigma)$. It can be shown that a necessary as well as
sufficient condition for torsion to be right-linear is that $\sigma$
satisfy the consistency condition
\be
\pi \circ (\sigma + 1) = 0.                                   \label{2.2.6}
\ee
Using the fact that $\pi$ is a projection one sees that the most general
solution to this equation is given by
\be
1 + \sigma = ( 1 - \pi) \circ \tau                            \label{2.2.7}
\ee
where $\tau$ is an arbitrary bilinear map
\be
\Omega^1(\c{A}) \otimes \Omega^1(\c{A})
\buildrel \tau \over \longrightarrow
\Omega^1(\c{A}) \otimes \Omega^1(\c{A}).                     \label{2.2.8}
\ee
If we choose $\tau = 2$ then we find $\sigma = 1 - 2 \pi$ and 
$\sigma^2 = 1$. The eigenvalues of $\sigma$ are then equal to $\pm 1$.  

This general formalism can be applied in particular to differential
calculi with a frame.  Since $\Omega^1(\c{A})$ is a free module the
maps $\sigma$ and $\tau$ can be defined by their action on the basis
elements:
\be
\sigma (\theta^a \otimes \theta^b) = 
S^{ab}{}_{cd} \theta^c \otimes \theta^d,   \qquad
\tau (\theta^a \otimes \theta^b) = 
T^{ab}{}_{cd} \theta^c \otimes \theta^d.                       \label{2.2.9}
\ee
By the sequence of identities
\be
f S^{ab}{}_{cd} \theta^c \otimes \theta^d = 
\sigma (f \theta^a \otimes \theta^b) = 
\sigma (\theta^a \otimes \theta^b f) =
S^{ab}{}_{cd} f \theta^c \otimes \theta^d                      \label{2.2.10}
\ee
and the corresponding ones for $T^{ab}{}_{cd}$ we conclude that the
coefficients $S^{ab}{}_{cd}$ and $T^{ab}{}_{cd}$ must lie in 
$\c{Z}(\c{A})$. From~(\ref{2.2.7}) the most general form for 
$S^{ab}{}_{cd}$ is
\be
S^{ab}{}_{cd} = 
T^{ab}{}_{ef} (\delta^e_c \delta^f_d - P^{ef}{}_{cd}) 
- \delta^a_c \delta^b_d.                                       \label{2.2.11}
\ee

A covariant derivative can be defined also by its action on the basis 
elements:
\be
D\theta^a = - \omega^a{}_{bc} \theta^b \otimes \theta^c.       \label{2.2.12}
\ee
The coefficients here are elements of the algebra. They are restricted
by (2.1) and the the two Leibniz rules. The torsion 2-form is defined as
usual as
\be
\Theta^a = d \theta^a - \pi \circ D \theta^a.             \label{deftorsion}
\ee
If $F^a{}_{bc} = 0$ then it is easy to check~\cite{DubMadMasMou96} 
that
\be
D_{(0)} \theta^a = - \theta \otimes \theta^a + 
\sigma (\theta^a \otimes \theta)                            \label{2.2.14}
\ee
defines a torsion-free covariant derivative.                    
The most general $D$ for fixed $\sigma$ is of the form
\be
D = D_{(0)} + \chi                                           \label{2.2.15}
\ee
where $\chi$ is an arbitrary bimodule morphism
\be
\Omega^1(\c{A}) \buildrel \chi \over \longrightarrow
\Omega^1(\c{A}) \otimes \Omega^1(\c{A}).               \label{2.2.16}
\ee
If we write
\be
\chi (\theta^a) = - \chi^a{}_{bc} \theta^b \otimes \theta^c   \label{2.2.17}
\ee
we conclude that $\chi^a{}_{bc} \in \c{Z}(\c{A})$. 
In general a covariant derivative is torsion-free provided the condition
\be
\omega^a{}_{de} P^{de}{}_{bc} = {1\over 2}C^a{}_{bc}          \label{2.2.19}
\ee
is satisfied. The covariant derivative~(\ref{2.2.15}) is torsion free 
if and only
if
\be
\pi \circ \chi = 0.                                            \label{2.2.20}
\ee

One can define a metric by the condition
\be
g(\theta^a \otimes \theta^b) = g^{ab}                           \label{2.2.21}
\ee
where the coefficients $g^{ab}$ are elements of $\c{A}$.  To be well
defined on all elements of the tensor product $\Omega^1(\c{A})
\otimes_\c{A} \Omega^1(\c{A})$ the metric must be bilinear and by
the sequence of identities
\be
f g^{ab} = g(f \theta^a \otimes \theta^b) 
= g(\theta^a \otimes \theta^b f) = g^{ab} f                     \label{2.2.22}
\ee
one concludes that the coefficients must lie in $\c{Z}(\c{A})$.
We define the metric to be symmetric if
\be 
g \circ \sigma \propto g.                                       \label{symm}
\ee
This is a natural generalization of the situation in ordinary
differential geometry where symmetry is respect to the flip which
defines the forms.  If $g^{ab} = g^{ba}$ then by a linear transformation
of the original $\lambda_a$ one can make $g^{ab}$ the components of the
Euclidean (or Minkowski) metric in dimension $n$. It will not
necessarily then be symmetric in the sense that we have just used the
word.

Introduce the standard notation $\sigma_{12} = \sigma \otimes 1$,
$\sigma_{23} = 1 \otimes \sigma$, to extend to three factors of a module
any operator $\sigma$ defined on a tensor product of two factors.
Then there is a natural continuation of the map (\ref{2.2.4}) to the 
tensor product $\Omega^1(\c{A}) \otimes_\c{A} \Omega^1(\c{A})$ 
given by the map
\be
D_2(\xi \otimes \eta) = D\xi \otimes \eta + 
\sigma_{12} (\xi \otimes D\eta).                                \label{2.2.4e}
\ee
The map $D_2 \circ D$ has no nice properties but if one introduces the
notation $\pi_{12} = \pi \otimes 1$ then by analogy with the commutative 
case one can set
\be
D^2 = \pi_{12} \circ D_2 \circ D
\ee
and formally define the curvature as the map
\be
\mbox{Curv} = D^2.                                              \label{curv}
\ee
Because of the condition (\ref{2.2.6}) Curv is left linear. It can
be written out in terms of the frame as
\be
\mbox{Curv} (\theta^a) =
- {1 \over 2} R^a{}_{bcd} \theta^c \theta^d \otimes \theta^b    \label{2.16}
\ee
Similarly one can define a Ricci map
\be
\mbox{Ric} (\theta^a) =
{1 \over 2} R^a{}_{bcd} \theta^c g(\theta^d \otimes \theta^b).  \label{2.17}
\ee
It is given by
\be
\mbox{Ric} \, (\theta^a) = R^a{}_b \theta^b.                   \label{2.18}
\ee
The above definition of curvature is not satisfactory in the
noncommutative case. For a discussion of this point we refer to the
article by Dubois-Violette {\it et al.}~\cite{DubMadMasMou96}.

The curvature of the covariant derivative $D_{(0)}$ defined in 
(\ref{2.2.14}) can be readily calculated. One finds after a short
calculation that it is given by the expression
\bea
\mbox{Curv} (\theta^a)\squeeze 
&& =\theta^2 \otimes \theta^a -
\pi_{12} \sigma_{23} \sigma_{12} 
(\theta^a \otimes \theta \otimes \theta) \nonumber\\[4pt]
&& = \theta^2 \otimes \theta^a +
\pi_{12} \sigma_{12}\sigma_{23} \sigma_{12}
(\theta^a \otimes \theta \otimes \theta).                      \label{2.19}
\eea
If  $\xi = \xi_a \theta^a$ is a general 1-form then since Curv is left
linear one can write
\bea
\mbox{Curv} (\xi ) \squeeze
&& = \xi_a \theta^2 \otimes \theta^a -
\pi_{12} \sigma_{23}  \sigma_{12}
(\xi \otimes \theta \otimes \theta) \nonumber\\[4pt]
&& = \xi_a \theta^2 \otimes \theta^a +
\pi_{12} \sigma_{12}\sigma_{23} \sigma_{12}
(\xi \otimes \theta \otimes \theta).                           \label{2.20}
\eea
We shall use this latter expression in Section~\ref{metrics}.

The compatibility of the covariant derivative~(\ref{2.2.12})  
with the metric is expressed by the condition~\cite{DubMadMasMou95}
\be
g_{23}\circ D_2= d\circ g .                                \label{met-comp}
\ee
The covariant derivative~(\ref{2.2.12}) is compatible with the metric if
and only if~\cite{DimMad96}
\be
\omega^a{}_{bc} + \omega_{cd}{}^e S^{ad}{}_{be} = 0.            \label{2.2.23}
\ee
This is a `twisted' form of the usual condition that
$g_{ad}\omega^d{}_{bc}$ be antisymmetric in the two indices $a$ and $c$
which in turn expresses the fact that for fixed $b$ the
$\omega^a{}_{bc}$ form a representation of the Lie algebra of the
Euclidean group $SO(N)$ (or the Lorentz group).  When $F^a{}_{bc} = 0$
the condition that~(\ref{2.2.12}) be metric compatible can be
written as
\be
S^{ae}{}_{df} g^{fg} S^{bc}{}_{eg} = g^{ab} \delta^c_d.         \label{2.2.24}
\ee

The algebra we shall consider is a $*$-algebra and we require that the
differential calculus be such that the reality condition 
\be
(df)^* = df^*                                                 \label{d-real}
\ee 
holds. Neither of the two differential calculi we shall introduce in
Section~4 satisfies this condition. In Section~7 we shall discuss how to
construct a real calculus $\Omega_R^*(\c{A})$ by taking a subalgebra 
of the tensor product of the two calculi. We shall require that for
arbitrary $f \in \c{A}$ and $\xi \in \Omega_R^1(\c{A})$ one has
\be
(f \xi)^* = \xi^* f^*, \qquad (\xi f)^* = f^* \xi^*.            \label{sola}
\ee
We shall suppose~\cite{DubMadMasMou95} that the extension of the
involution to the tensor product is given by
\be
(\xi \otimes \eta)^* = \sigma(\eta^* \otimes \xi^*).             \label{TPI}
\ee
A change in $\sigma$ therefore implies a change in the definition of
an hermitian tensor. The reality condition for the metric becomes 
then~\cite{FioMad98}
\be
g\circ \sigma (\eta^* \otimes \xi^*) = 
(g(\xi \otimes \eta))^*                                      \label{reality}
\ee 
There is also a reality condition on the covariant derivative and the
curvature~\cite{FioMad98} which imply that the generalized flip 
$\sigma$ must satisfy the braid equation.

\sect{The quantum Euclidean space}

The $R$-matrix or braid matrix $\hat R\equiv\Vert \hat R^{ij}_{hk}\Vert$ 
of $SO_q(3)$ \cite{frt} is a $3^2\times 3^2$ matrix satisfying
the braid equation
\be
\hat R_{12}\,\hat R_{23}\,\hat R_{12} =\hat R_{23}\,\hat R_{12}\,
\hat R_{23}.                                                  \label{braid1}
\ee
Here we have used the conventional tensor notation
$\hat R_{12} = \hat R \otimes 1$ used above for $\sigma$.
By repeated application of the Equation~(\ref{braid1}) one finds 
\be
f(\hat R_{12})\,\hat R_{23}\,\hat R_{12} 
=\hat R_{23}\,\hat R_{12}\,f(\hat R_{23})                     \label{braid2}
\ee
for any polynomial function $f(t)$ in one variable.  The
Equations~(\ref{braid1}) and~(\ref{braid2}) are evidently satisfied also
after the replacement $\hat R\rightarrow\hat R^{-1}$.  The braid matrix
admits the projector decomposition
\be
\hat R = q\c{P}_s - q^{-1}\c{P}_a + q^{-2}\c{P}_t             \label{decompo}
\ee
with 
\be
\c{P}_{\mu}\c{P}_{\nu} = \c{P}_\mu \delta_{\mu\nu}, \qquad
\sum_{\mu}\c{P}_{\mu} = 1, \quad \mu,\nu = s,t,a.
\ee
The $\c{P}_s$, $\c{P}_a$, $\c{P}_t$ are $SO_q(3)$-covariant
$q$-deformations of respectively the symmetric trace-free,
antisymmetric and trace projectors. The trace projector is
1-dimensional and its matrix elements can be written in the form
\be
\c{P}_t{}_{kl}^{ij} = (g^{mn}g_{mn})^{-1} g^{ij}g_{kl},         \label{bibi}
\ee
where~\cite{frt} the `metric matrix' $g= (g_{ij})$ is a
$SO_q(3)$-isotropic tensor, a deformation of the metric matrix on the
classical Euclidean space.  The $\hat R$ and $g$ satisfy the relations
\be
g_{il}\,\hat R^{\pm 1}{}^{lh}_{jk} = 
\hat R^{\mp 1}{}^{hl}_{ij}\,g_{lk}, \qquad
g^{il}\,\hat R^{\pm 1}{}_{lh}^{jk} = 
\hat R^{\mp 1}{}_{hl}^{ij}\,g^{lk}.                           \label{utile1}
\ee
The lower-case Latin indices $i,j,\dots$ will take
the values $(-, 0, +)$.  The quantum Euclidean space is the formal
(associative) algebra $\c{A}$ with generators $x^i = (x^-, x^0, x^+)$ 
and relations
\be
\c{P}_a{}_{kl}^{ij} x^k x^l=0                                 \label{xxcr}
\ee
for all $i,j$. We introduce a grading in $\c{A}$ by requiring that the
degree of $(x^-, x^0, x^+)$ be respectively equal to $(-1, 0, +1)$. The
matrix elements of $\hat R$ and therefore of all the projectors fulfill
the condition
\be
\hat R^{ij}_{kl} = 0 \qquad \mbox{if} \quad i+j \neq k+l.
\ee
Consequently, all the terms appearing on the left-hand side of 
Equation~(\ref{xxcr}) have the same total degree $i\!+\!j$. 
If we use the explicit expression~\cite{frt} for $\c{P}_a$ and set
\be
h = \sqrt q - 1/\sqrt q, 
\ee
then the relations~(\ref{xxcr}) can be written in the form
\be
\ba{l}
x^- x^0 = q\, x^0 x^-,\\[2pt]
x^+ x^0 = q^{-1} x^0 x^+,\\[2pt]
[x^+, x^-] = h (x^0)^2.
\ea                                                        \label{xxcr3}
\ee
The first two equations define two copies of the $q$-quantum plane with
$q$ and $q^{-1}$ as deformation parameter and a common generator $x^0$.

The metric matrix is given by 
$g_{ij} = g^{ij}$ with
\be
g_{ij} = \left(
\ba{ccc}
0            &0 & 1/\sqrt{q}  \\
0            &1 &    0        \\
\sqrt{q}     &0 &    0        \\
\ea\right).                                                 \label{metric}
\ee
When $q \in \b{R}^+$ one obtains the real quantum
Euclidean space by giving $\c{A}$ the structure of a $*$-algebra with
\be
(x^-)^* = \sqrt q x^+, \qquad (x^0)^* = x^0, \qquad
(x^+)^* = {1\over\sqrt q} x^-.
\ee
This can be written in terms of the metric as
\be
(x^i)^* = x^j g_{ji}.                                         \label{real}
\ee
We can use the summation convention if we consider the involution to
lower (or raise) an index.  The condition (\ref{real}) is an
$SO_q(3,\b{R})$-covariant condition and three linearly independent,
hermitian `coordinates' can be obtained as combinations of the
$x^i$. We define
\be
x^r = \Lambda^r_i x^i,
\qquad \Lambda^r_i:= {1 \over \sqrt2} \left(
\ba{ccc}
    1    &0    &\sqrt q\\
    0 &\sqrt  2   &0\\
    i    &0    &-i\sqrt q
\ea\right).                                                  \label{realmat}
\ee
With respect to the new `coordinates' the metric is given by
$$
g^{rs} = g^{ij} \Lambda^r_i \Lambda^s_j = {1 \over 2} \left(
\ba{ccc}
   q+1   &0  &i(q-1)\\
    0    &2  &0\\
 -i(q-1) &0  &q+1
\ea\right).
$$
The metric is hermitian but no longer real. In the limit $q\to 1$ one
sees that $g^{rs} \to \delta^{rs}$.  It is more convenient to remain
with the original `coordinates' and a real metric.

The `length' squared
\be
r^2 := g_{ij} x^i x^j 
\ee
is an $SO_q(3,\b{R})$-invariant, real and generates the center
$\c{Z}(\c{A})$ of $\c{A}$. Using~(\ref{real}) and (\ref{xxcr3}) it can
be written also in the forms
\be
r^2 = (x^i)^* x^i = (\sqrt{q}+\frac 1{\sqrt{q}})x^-x^++q 
(x^0)^2 = (\sqrt{q}+\frac 1{\sqrt{q}})x^+x^-+q^{-1} (x^0)^2.
\label{squarelenght}
\ee
The commutation relation between $x^+$ and
$x^-$ can also be written in the form
\be
q x^+ x^- - q^{-1} x^- x^+ = h r^2.
\ee
There is obviously no obstruction to extending $\c{A}$ by adding to it
the inverse $r^{-2}$ of $r^2$, and also the square roots $r, r^{-1}$ of
these two elements. We shall further extend the algebra $\c{A}$ by
adding the inverse $(x^0)^{-1}$ of $x^0$ as a new generator with the
obvious extra relations.

Finally, we observe that if $|q|=1$ a compatible involution is defined
rather by $(x^i)^*= x^i$. The algebra describes a quantum space
covariant under the noncompact section $SO_q(2,1)$ of $SO_q(3,\b{C})$.

\sect{The $q$-deformed Euclidean differential calculi}        \label{calculi}

There are two $SO_q(3)$-covariant differential calculi~\cite{cawa}
$\Omega^*(\c{A})$ and $\bar \Omega^*(\c{A})$ over $\c{A}$ neither of
which satisfy the reality condition~(\ref{d-real}). In this section we
shall write them in a form which will permit us in Section~7 to look
for a sensible definition of a real differential calculus. Let $d$ and
$\bar d$ be the respective differentials and set $\xi^i = dx^i$ and
$\bar \xi^i = \bar dx^i$. The calculi are determined respectively by
the commutation relations
\be
x^i \xi^j = q\,\hat R^{ij}_{kl} \xi^k x^l                     \label{xxicr}
\ee
for $\Omega^1(\c{A})$ and
\be
x^i \bar \xi^j = 
q^{-1}\,\hat R^{-1}{}^{ij}_{kl} \bar\xi^k x^l                  \label{xxicr'}
\ee
for $\bar\Omega^1(\c{A})$.  Using formulae~(\ref{utile1}) and the fact
that we would like to require that the calculus satisfy~(\ref{d-real}),
it is easy to see~\cite{olezu} that, if $q\in \b{R}^+$, 
it is not possible to extend the
involution~(\ref{real}) to either calculus. This is possible only if
$|q|=1$, and will be considered elsewhere. On the other hand, if 
$q\in \b{R}^+$ the involution can be extended to the direct sum
$\Omega^1(\c{A})\oplus\bar\Omega^1(\c{A})$ by setting, as well as
(\ref{sola}),
\be
(\xi^i)^* = \bar\xi^j g_{ji},                               \label{barstern}
\ee
since this exchanges relations~(\ref{xxicr}) and~(\ref{xxicr'}).  
This will be equivalent to
\be
(df)^* = \bar d f^*.                                     \label{real-cond}
\ee
If in the limit $q\to 1$ we identify $\bar\xi^i=\xi^i$,
we recover the standard involution and the
differential is real.

By means of Equations~(\ref{xxicr}) and~(\ref{utile1}) it is 
straightforward to check that the commutation relations
\be
r^2 \xi^i = q^2\,\xi^i r^2.                                    \label{balla}
\ee
are satisfied.  The commutation relations between the $\xi^i$ are
derived by taking the differential of~(\ref{xxicr}):
\be
\c{P}_s{}_{kl}^{ij} \xi^k \xi^l = 0, \qquad\qquad
\c{P}_t{}_{kl}^{ij} \xi^k \xi^l = 0.                          \label{xixicr}
\ee
If we extend the grading of $\c{A}$ to $\Omega^1(\c{A})$ by setting
$\mbox{deg}(\xi^i)= \mbox{deg}(x^i)$, we find that each term in
(\ref{xxicr}), (\ref{xixicr}) must have the same total degree. 

Written out explicitly the Equations~(\ref{xxicr}) become
\bea
&&x^- \xi^- = q^2\, \xi^- x^-, \nonumber\\[2pt]
&&x^0 \xi^- = q\, \xi^- x^0, \\[2pt]
&&x^+ \xi^- = \xi^- x^+,   \nonumber\\[5pt]
&&x^- \xi^0 = q\,\xi^0 x^- +(q^2-1)\xi^-x^0, \nonumber\\[2pt]
&&x^0 \xi^0 = q\,\xi^0 x^0 - h (q+1) \xi^- x^+, \\[2pt]
&&x^+ \xi^0 = q\,\xi^0 x^+, \nonumber\\[5pt]
&&x^- \xi^+ = \xi^+ x^- - h (q+1) \xi^0 x^0 + 
h^2 (q+1) \xi^- x^+, \nonumber\\[2pt]
&&x^0 \xi^+ = q\, \xi^+ x^0 + (q^2-1)\xi^0 x^+,\\[2pt]
&&x^+ \xi^+ = q^2\, \xi^+ x^+ \nonumber
\eea
and the Equations~(\ref{xixicr}) become~\cite{fiodet} 
\bea
&&(\xi^{\pm})^2 = 0,\nonumber\\[2pt]
&&(\xi^0)^2 = h \xi^- \xi^+,\nonumber\\[2pt]
&&\xi^- \xi^0 = - q^{-1}\xi^0 \xi^-,\\[2pt]
&&\xi^0 \xi^+ = - q^{-1}\xi^+ \xi^0,\nonumber\\[2pt]
&&\xi^- \xi^+ = - \xi^+ \xi^-.\nonumber
\eea

Consider the $SO_q(3)$-invariant 1-form
\be
\eta := g_{ij} x^i \xi^j = q^{-1}\, g_{ij} \xi^j x^i.
\ee
Using the projector decomposition of the $\hat R$-matrix, the
relations~(\ref{utile1}) between the $\hat R$-matrix and its inverse
as well as the relations~(\ref{xxcr}) which define the algebra one can
easily verify that
\be
[\eta, x^i] = - q^{-2} (q-1) r^2\xi^i.
\ee
Hence we conclude that
\be
\theta := (q-1)^{-1} q^2 r^{-2} \eta                        \label{deftheta}
\ee
is the `Dirac operator'~(\ref{dirac}) of $\Omega^1(\c{A})$.  It
satisfies the conditions
\be
d\theta = 0, \quad \theta^2 = 0.                               \label{fifa}
\ee
It is of interest to note that 
\be
dr^2 = (1-q^{-2})\, r^2\,\theta.                              \label{diff-r}
\ee
Since $r \in \c{Z}(\c{A})$, from relations (\ref{diff-r}) one
concludes immediately that the differential calculus cannot be based
on derivations as outlined in the first section. That is, there can
exist no decompositions of $\theta$ as in~(\ref{dirac}). This fact is
not necessarily a defect but we shall change it anyway by adding an
extra element, the `dilatator', to the algebra since at the same time
we can reduce $\c{Z}(\c{A})$ to $\b{C}$.

The most general Ansatz for the frame can be written in the form
$\theta^a := \theta^a_i \xi^i$ where the $\theta^a_i$ are elements of
$\c{A}$. We shall let the lower-case Latin indices $a,b,...$ take the
same values $(+,0,-)$ as $i,j,...$. The condition~(\ref{fund}) implies
that
\[
(r^2\theta^a_i - q^{-2} \theta^a_ir^2) \xi^i = 0
\]
from which we can conclude that
\[
r^2\theta^a_i - q^{-2} \theta^a_i r^2= 0.
\]
This equation has obviously no solution since $r^2 \in \c{Z}(\c{A})$. 
To remedy this problem we extend the algebra $\c{A}$ by adding an
extra generator $\Lambda$, the `dilatator', and its inverse
$\Lambda^{-1}$, chosen such that
\be
x^i \Lambda = q\, \Lambda x^i .                         \label{lambda}
\ee
The introduction of a new generator $\Lambda$ is necessary also in a
different context, namely in the inhomogeneous extension of the
homogeneous quantum groups $SL_q(N)$, $SO_q(N)$ and
$q$-Lorentz~\cite{wessvari, munich, maj2}; more precisely, $\Lambda$
appears in the coproduct of the translation part generators.  We do
not know if this is a coincidence or there is some more fundamental
link between the two phenomena.

We extend the original algebra $\c{A}$ defined by 
Equations~(\ref{xxcr3}) by the addition not only of 
$(x^0)^{-1}, r^{\pm 1}$ but also of $\Lambda^{\pm 1}$ as new 
generators.
Since
$$
r \Lambda = q\, \Lambda r,
$$
clearly the center of $\c{A}$ is now trivial: $\c{Z}(\c{A}) = \b{C}$.
It is natural to extend the grading by setting $\mbox{deg}(\Lambda) = 0$. 
 We shall choose $\Lambda$ to be
unitary
\be
\Lambda^* = \Lambda^{-1}.
\ee
To within a normalization this is the only choice compatible
with the commutation relations.

We can now consider the previous Ansatz for the frame but with
coefficients in the algebra $\c{A}$. We must assume a linear
dependence on the generator $\Lambda$ and write
\be 
\theta^a := \Lambda^{-1}\, \theta^a_j \xi^j                      \label{ok}
\ee
where the $\theta^a_j$ are elements of $\c{A}$ which do not depend on
$\Lambda$. The $\Lambda$-dependence is here dictated by the condition 
$[r,\theta^a] = 0$. The condition $[x^i, \theta^a] = 0$ becomes
\be
x^i \theta^a_j = \hat R^{-1}{}^{ki}_{lj}\theta^a_k x^l.          \label{xdcr}
\ee
Written out explicitly these Equations become
\bea
&&x^- \theta^a_- = q^{-1} \theta^a_- x^- - 
q^{-1} (q^2-1) \theta^a_0 x^0 + h^2 (1+q) \theta^a_+ x^+,\nonumber\\[2pt]
&&x^0 \theta^a_- = \theta^a_- x^0 + h(1+q) \theta^a_0 x^+, \label{1st}\\[2pt]
&&x^+ \theta^a_- = q\,\theta^a_- x^+,\nonumber\\[5pt]
&&x^- \theta^a_0 = \theta^a_0 x^- + h(1+q) \theta^a_+ x^0,\nonumber\\[2pt]
&&x^0\theta^a_0 = \theta^a_0 x^0 - 
q^{-1} (q^2-1) \theta^a_+ x^+,                        \label{2nd}\\[2pt]
&&x^+ \theta^a_0 = \theta^a_0 x^+,\nonumber\\[5pt]
&&x^- \theta^a_+ = q\, \theta^a_+ x^-,\nonumber\\[2pt]
&&x^0 \theta^a_+ = \theta^a_+ x^0,                      \label{3rd}\\[2pt]
&&x^+ \theta^a_+ = q^{-1} \theta^a_+ x^+.\nonumber
\eea

The Equations~(\ref{3rd}) contain only $\theta^a_+$ and admit, to within
an arbitrary factor in the center of $\c{A}$, only two independent
solutions. We choose
\beas
&&\theta^-_+ = 0,                                  \\[2pt]
&&\theta^0_+ = 0,                                  \\[2pt]
&&\theta^+_+ = r^{-2} x^0,
\eeas
so that $\theta^+$ contains a term proportional to $\xi^+$. 
We replace this result in Equations~(\ref{2nd})
which become then equations for $\theta^a_0$ and we find the solutions
\beas
&&\theta^-_0 = 0,                                    \\[2pt]
&&\theta^0_0 = r^{-1},                               \\[2pt]
&&\theta^+_0 = - (q+1)r^{-2} x^+.
\eeas
Finally, the Equations~(\ref{1st}) can be solved for $\theta^a_-$
yielding
\beas
&&\theta^-_- = (x^0)^{-1},                                     \\[2pt]
&&\theta^0_- = \sqrt q\, (q+1) r^{-1} (x^0)^{-1} x^+,          \\[2pt]
&&\theta^+_- = - \sqrt q\, q (q+1) r^{-2} (x^0)^{-1} (x^+)^2.
\eeas
If we place these coefficients into Equation~(\ref{ok}) we completely
determine the frame.
We shall include a rescaling by a normalization factor
$\alpha$, which we leave free for the moment:
\bea
&&\theta^- =  \alpha^{-1}\Lambda^{-1}\, (x^0)^{-1} \xi^-,  \nonumber\\[4pt]
&&\theta^0 =  \alpha^{-1}\Lambda^{-1}\, r^{-1}(\sqrt q\, 
(q+1) (x^0)^{-1} x^+ \xi^- + \xi^0), \label{stehbein3}\\[4pt]
&&\theta^+ =  \alpha^{-1}\Lambda^{-1}\, r^{-2}(-\sqrt q\, 
q (q+1) (x^0)^{-1} (x^+)^2 \xi^- - (q+1) x^+ \xi^0 + x^0 \xi^+).\nonumber
\eea
Note that we have enumerated the $\theta^a$ so that
$\mbox{deg}(\theta^a)=(-1,0,1)$, exactly as for the $x^i$ and $\xi^i$.

The above $\theta^a$ are determined up to linear transformations 
with coefficients depending on $r$. 
The commutation rules between $\Lambda$ and $\theta^a$
will depend on their $r$-normalization
as well as the commutation rules between $\Lambda$ 
and $\xi^i$, which we have not specified yet. 
By differentiating Equation~(\ref{lambda}) we obtain the condition
\be
\xi^i \Lambda + x^i d\Lambda = q d\Lambda x^i + q \Lambda \xi^i
\ee
on the differential $d\Lambda$. A possible solution is given by the two
conditions 
\be
x^i\, d\Lambda = q d\Lambda x^i, \qquad
\xi^i \Lambda = q \Lambda \xi^i.                                \label{zero}
\ee
In particular one can consistently set $d\Lambda = 0$. We shall do in
the sequel although it means that the condition $df = 0$ does not imply
that $f \in \b{C}$. This is not entirely satisfactory since one
would like the only `functions' with vanishing exterior derivative to be
the `constant functions'. It could be remedied by considering
a more general solution to Equation~(\ref{zero}). This would however
complicate our calculations since it would increase the number of
independent forms by one.  A necessary condition for $d\Lambda=0$ is
that $[\Lambda, \lambda_a] = 0$. The $\lambda_a$ which we give below in
Equation~(\ref{lambda'}) satisfy the latter because they are homogeneous
functions of $x^i$ of degree zero. The condition that $\Lambda$ commutes
with the $\theta^a$ fixes the normalization factor $\alpha$ in
(\ref{stehbein3}) to be proportional to 1.  In Section~\ref{climit} we
shall choose $\alpha$ in $\b{R}^+$, first as a fixed scale and then as a
suitable function of $h$, so as to adjust the commutative limit of the
frame.  Another possibility would be to impose a commutation relation
between $\xi$ and $\Lambda$ of the type
\be
\Lambda \xi^i= \xi^i\Lambda.                                \label{zero'}
\ee
This follows, for example, from the condition~\cite{oleg2} 
$\Lambda df = q d \Lambda f$ which in turn would be equivalent to
considering $\Lambda$ a suitable element of the associated Heisenberg
algebra\footnote{Such an element has a simple realization in terms of
$x^i$ and $SO_q(3)$-covariant twisted derivatives~\cite{oleg2}}.  The
condition that $\Lambda$ commutes with the $\theta^a$ would now fixes
the normalization factor $\alpha$ in (\ref{stehbein3}) to be
proportional to $r^{-1}$ and will thus change the metric by a
conformal factor $r^2$.  After imposing either commutation rule
between $\Lambda$ and $\theta^a$ the frame will be determined up to a
linear transformation with coefficients in $\b{C}$.
%%%%%%%%%%%%%
%
% The condition $\Lambda df = q d \Lambda f$ is really odd! Cf. the
% commutative limit...
%
%%%%%%%%%%%%%
A direct and lengthy calculation\footnote{Performed using the program
for symbolic computations REDUCE.} with the explicit expression for the
matrix $\hat R^{ab}_{ij}$~\cite{frt} shows that the $\theta^a_i$ satisfy
the relations
\be
\hat R^{ab}_{cd}\, \theta^d_j\theta^c_i = 
\theta^b_l\theta^a_k\, \hat R^{kl}_{ij}.                   \label{basic}
\ee
These are $3^4=81$ equations. However, since both 
$R^{kl}_{ij}\equiv \hat R^{lk}_{ij}$ and $\theta^a_i$ are
lower-triangular matrices, $3^2(3^2-1)/2 = 36$ of them are trivial
identities. The proof actually consists then in checking `only' 45
equations.  By repeated application of relations (\ref{basic}) it
immediately follows that the same relations hold also after replacing
$\hat R$ by any polynomial $f(\hat R)$,
\be
f(\hat R)^{ab}_{cd}\, \theta^d_j\theta^c_i = 
\theta^b_l\theta^a_k\, f(\hat R)^{kl}_{ij}.                    \label{basic1}
\ee
In particular we can choose $f(\hat R) = \c{P}_t,\c{P}_s$.  With the
help of relations~(\ref{fund}) and~(\ref{ok}) we find as a consequence
that the $\theta^a$ have the same commutation relations as the $\xi^i$:
\be
\c{P}_t{}^{ab}_{cd} \theta^c \theta^d = 0 \qquad\qquad
\c{P}_s{}^{ab}_{cd} \theta^c \theta^d = 0.                 \label{ththcr}
\ee
Therefore the $P$ of Equation~(\ref{manca}) is given by
\be
P = \c{P}_a.                                               \label{whynot}
\ee
It is the $SO_q(3)$-covariant deformation of the antisymmetric projector.  
Note also that for $f(\hat R)=\c{P}_t$ Equation~(\ref{basic1})
is equivalent to the relations
\be
g_{cd}\theta^d_j\theta^c_i = \kappa g_{ij}\qquad\qquad
g^{ij}\theta^b_j\theta^a_i = \kappa g^{ab}               \label{mamam}
\ee
with 
$$
\kappa = (g_{kl}g^{kl})^{-1} g_{ab} \theta^b_i \theta^a_j g^{ji}
= r^{-2}\alpha^{-2}
$$
because of Equation~(\ref{bibi}).

Consider the elements $\lambda_a \in \c{A}$ with
\bea
&&\lambda_- =   h^{-1}q \Lambda(x^0)^{-1}\alpha  x^+,      \nonumber\\[2pt]
&&\lambda_0 = 
- h^{-1}\sqrt q \Lambda\alpha  (x^0)^{-1} r,         \label{lambda'}\\[2pt]
&&\lambda_+ = - h^{-1}\Lambda\alpha  (x^0)^{-1} x^-.       \nonumber
\eea
By direct calculation one verifies that 
$\theta = - \lambda_a \theta^a$ is given by~(\ref{deftheta}).
Since $\Lambda$ is unitary the hermitian adjoints $\lambda_a^*$ 
are given by
\be
\lambda_\pm^* = - \Lambda^{-2} g^{\pm b}\lambda_b, \qquad
\lambda_0^* = \Lambda^{-2} g^{0b}\lambda_b.                    \label{stira}
\ee
The fact that the $\lambda_a$ are not anti-hermitian is related to the
fact that the differential $d$ is not real.  We have chosen this rather
odd normalization to have the commutation relations (\ref{lambda-com})
below.  A straightforward calculation yields the commutation relations
$[\lambda_a, x^i] = q\,\Lambda e_a^i$ with
\be
(e_a^i) = \,\alpha\,\left(
\ba{ccc}
+ x^0 & 0 & 0 \\
- (\sqrt q + 1/\sqrt q) x^+  & r & 0\\
-\sqrt q\, (q+1) (x^0)^{-1} (x^+)^2 & (q+1) r (x^0)^{-1} x^+ 
& r^2 (x^0)^{-1}
\ea\right).
\ee
The $e_a^i$ is the inverse matrix of $\theta^a_i$:
\be
e_a^i\theta^a_j = \delta^i_j\qquad\qquad 
\theta^a_je^j_b=\delta^a_b.                                 \label{inverse}
\ee
Relations (\ref{basic}), (\ref{mamam}) imply that the matrix elements
$e^i_a$ satisfy also the `$RTT$-relations' \cite{frt}
\be
\hat R_{kl}^{ij}\, e^k_a e^l_b
= e^i_c e^j_d \,\hat R_{ab}^{cd}                                \label{RTT}
\ee
as well as the `$gTT$-relations' \cite{frt}
\be
g^{ab} e^i_a e^j_b = r^2 \alpha^2 g^{ij}\qquad\qquad
g_{ij} e^i_a e^j_b = r^2 \alpha^2 g_{ab}.                        \label{gTT}
\ee
Thus $(\alpha r)^{-1}e^i_a$ fulfill the same commutation relations of
the generators $T^i_a$ of $SO_q(3)$ and in this sense may be seen as a
`local' realization of $T^i_a$. However they do not satisfy the same
$*$-relations, nor is there an analog for the coproduct of $SO_q(3)$.
The $\lambda_a$ satisfy the commutation relations
\bea
&&\lambda_- \lambda_0 = q\, \lambda_0 \lambda_-,   \nonumber\\[2pt]
&&\lambda_+ \lambda_0 = q^{-1} \lambda_0 \lambda_+,\label{lambda-com}\\[2pt]
&&[\lambda_+, \lambda_-] = h\, (\lambda_0)^2.      \nonumber
\eea
These are the same commutation relations as those satisfied
by the $x^i$. This is a remarkable fact and it underlines how weak a
constraint the commutation relations are on the algebra. The
$\lambda_a$ are related in fact to the $x^i$ by a rather complicated
nonlinear relation (\ref{lambda'}).  They differ however from the
$x^a$ in that they commute with $\Lambda$.

The Equations~(\ref{lambda-com}) can be rewritten more compactly in the
form
\be
\c{P}_a{}^{ab}_{cd} \lambda_a\lambda_b = 0,               \label{lambdacr1}
\ee
so that the constants $F^a{}_{bc}$ and $K_{ab}$ in Equation~(\ref{manca}) 
vanish. Hence Equation~(\ref{rot-coef}) is satisfied with
\be
C^a{}_{bc} = - 2 \lambda_d \c{P}_a{}^{(da)}{}_{bc}
\ee
an expression which is consistent with (\ref{consis1}) because of 
(\ref{whynot}). Finally, it is easy to check that
\be
g^{ab} \lambda_a \lambda_b = q h^{-2} (\Lambda\alpha)^2.    \label{lambdacr2}
\ee

One can easily find the relation between the $SO_q(3)$-covariant
derivatives $\partial_i$ introduced in Ref. \cite{cawa}, which fulfil the
modified Leibniz rule
\be
\partial_ix^j=\delta_i^j+q\hat R^{jh}_{ik}x^k\partial_h ,    
\ee
and the $e_a$.
From the decomposition $d=\xi^i\partial_i=\theta^ae_a$ it is evident that
\be
\partial_i= \Lambda^{-1}\tilde\theta^a_ie_a,              \label{covder}
\ee
where we have now decompesed $\theta^a$ in the form
$\theta^a=\xi^i\Lambda^{-1}\tilde\theta^a_i$.

We conclude this section by listing the basic formulae which, besides
(\ref{xxicr'}), characterize the barred differential calculus
$\bar\Omega^1(\c{A})$.  The analogs of (\ref{xixicr}) are obtained by
taking the differential of~(\ref{xxicr'}):
\be
\c{P}_s{}_{kl}^{ij} \bar\xi^k \bar\xi^l = 0, \qquad\qquad
\c{P}_t{}_{kl}^{ij} \bar\xi^k \bar\xi^l = 0.             \label{barxibarxicr}
\ee
All relations are compatible with the grading of $\c{A}$ extended by
setting $\mbox{deg}(\bar\xi^i)=\mbox{deg}(x^i)$.  The
$SO_q(3)$-invariant 1-form
\be
\bar\eta := g_{ij} x^i \bar\xi^j = q\,g_{ij} \bar\xi^j x^i   
\ee
fulfills in $\bar\Omega^1(\c{A})$  the commutation relations
\be
[\bar\eta, x^i] = - q^2 (q^{-1}-1) r^2 \bar\xi^i.
\ee
Hence
\be
\bar\theta := (q^{-1} -1)^{-1} q^{-2}  r^{-2} \bar\eta    \label{defbartheta}
\ee
is the `Dirac operator'~(\ref{dirac}) of $\bar\Omega^1(\c{A})$, and
\be
\bar d\bar\theta = 0, \quad \bar\theta^2 = 0.
\ee
 From the definitions of $\theta$, $\bar\theta$ and the involution
it follows that in $\Omega^1(\c{A})\oplus\bar\Omega^1(\c{A})$
\be
\theta^* = - \bar\theta, \qquad (\bar\theta)^* = - \theta.     \label{stern}
\ee
It is straightforward to show that
\be
\bar dr^2 = (1-q^2)\, r^2\,\bar\theta.
\ee
There exist a frame $\bar\theta^a$ for the differential calculus 
$\bar\Omega^*(\c{A})$ in the form
\be 
\bar\theta^a := \Lambda  \bar\theta^a_j\bar\xi^i,        \label{barok}
\ee
where the $\bar\theta^a_j$ are elements of $\c{A}$ which do not
depend on $\Lambda$. The $\Lambda$-dependence is here again dictated
by the condition $[r,\bar \theta^a] = 0$. The condition
$[x^i, \bar \theta^a] = 0$ becomes
\be
\bar \theta^a_i x^j = 
\hat R^{jk}_{il} x^l \bar \theta^a_k.                 \label{barxdcr}
\ee
The elements $\bar\lambda_a\in\c{A}$ are introduced
through the decomposition
\be
\bar\theta = - \bar\lambda_a \bar\theta^a                  \label{bardirac}
\ee
and are dual to $\bar\theta^a$.  As above, $\bar \theta^a$ and the
corresponding $\bar\lambda_a$ are determined up to a linear
transformation with coefficients in $\b{C}$. Now note that from $0 =
[\c{A},\theta^a]^* = [(\theta^a)^*,\c{A}]$ it follows that for
$q$ real positive 
$(\theta^a)^*$ is a combination of $\bar\theta^b$. We choose the
second basis so that
\be
(\theta^a)^*= \bar\theta^b g_{ba}.                          \label{bar-theta}
\ee
This will automatically yield $\mbox{deg}(\bar\theta^a)=(-1,0,1)$ 
and the same commutation relations as for the $\theta^a$,
because of relations (\ref{utile1}):
\be
\c{P}_t{}^{ab}_{cd} \bar\theta^c \bar\theta^d = 0 \qquad\qquad
\c{P}_s{}^{ab}_{cd} \bar\theta^c \bar\theta^d = 0.          \label{bth-bthcr}
\ee
The explicit expressions for $\bar \theta^a$ becomes
\bea
&&\bar\theta^- = (\alpha q)^{-1} \, \Lambda \, r^{-2}(x^0 \bar\xi^- - 
(q^{-1}+1) x^- \bar\xi^0 -\frac1{\sqrt q}\, q^{-1}(q^{-1}+1) 
(x^0)^{-1} (x^-)^2 \bar\xi^+),                        \nonumber\\[4pt]
&&\bar\theta^0 = (\alpha q)^{-1}\, \Lambda \, r^{-1} (\bar\xi^0 + 
\frac 1{\sqrt q}\, 
(q^{-1}+1)(x^0)^{-1} x^- \bar\xi^+),        \label{barstehbein3}\\[4pt]
&&\bar\theta^+ = (\alpha q)^{-1}\, \Lambda \, (x^0)^{-1} \bar \xi^+.   
\nonumber   
\eea
The corresponding $\bar\lambda_a$ are given by 
$$
\bar\lambda_\pm = \Lambda^{-2} \lambda_\pm, \qquad
\bar\lambda_0 = - \Lambda^{-2} \lambda_0
$$
and therefore the involution on the $\lambda_a$ becomes
\be
\lambda_a^* = - g^{ab}\bar\lambda_b.
\ee
This is to be compared with~(\ref{real}).

\sect{Metrics and linear connections on the \\ 
quantum Euclidean space}                                     \label{metrics}

We now look for metrics, generalized permutations and covariant
derivatives corresponding to each of the above differential
calculi. They will be essential ingredients do determine the correct
correspondence between mathematical objects and physical observables.
We shall see that it is not possible to satisfy all the requirements of
Section~\ref{preli}. Nevertheless, leaving the problem of reality aside
for the moment, to be treated in Section~\ref{involu}, we show that if
we allow a conformal factor in the metric then a unique linear
connection exists which is metric compatible and automatically
$SO_q(3)$-covariant.

We define covariant derivatives $D$ on $\Omega^1(\c{A})$ and $\bar D$ on
$\bar\Omega^1(\c{A})$ as maps
\beas
&&\Omega^1(\c{A}) \buildrel D \over \longrightarrow
\Omega^1(\c{A}) \otimes \Omega^1(\c{A}), \\[2pt]
&&\bar\Omega^1(\c{A}) \buildrel \bar D \over \longrightarrow
\bar\Omega^1(\c{A}) \otimes \bar\Omega^1(\c{A})
\eeas
which satisfy left and right Leibniz rules. 

In accord with Equation~(\ref{2.2.9}), we look
for a generalized permutation $\sigma$ by starting with the Ansatz
\be
\sigma (\theta^a\otimes \theta^b) =
S^{ab}{}_{cd}\,\theta^c\otimes \theta^d                     \label{gigio2}
\ee
with $S^{ab}{}_{cd}$ complex numbers and we impose the
condition (\ref{2.2.6}). The Equation (\ref{ththcr})
implies that 
\be
S = C_s\c{P}_s - \c{P}_a + C_t\c{P}_t,                        \label{gala}
\ee
where $C_s$ and $C_t$ are complex $3^2\times 3^2$ matrices.
Similarly, according to (\ref{2.2.21}) to define a
metric $g$ we start with the Ansatz
\be
g(\theta^a\otimes\theta^b) = g^{a b}                         \label{fine`}
\ee
with $g^{ab}$ again complex numbers.  In view of Equation~(\ref{utile1})
a necessary condition for the metric-compatibility
condition~(\ref{2.2.24}) is that $g_{ab}$ be proportional to the matrix
defined in (\ref{metric}) and either $S = \hat R$ or $S = \hat R^{-1}$.
Without loss of generality we can set the proportionality factor equal
to 1, since this amounts to a redefinition of the factor $\alpha$ in
(\ref{stehbein3}).  But because of (\ref{decompo}) these forms for $S$
are clearly not compatible with (\ref{gala}).

The best one can do is to weaken (\ref{2.2.24}) to a condition of
proportionality.  To fulfill the latter it is sufficient that $S$ be
proportional to $\hat R$ or $\hat R^{-1}$.  Then the double requirement
admits the two solutions
\be
S = q\hat R, \qquad\qquad S = (q\hat R)^{-1},                 \label{double}
\ee
corresponding respectively to $C_s = q^2$, $C_t=q^{-1}$ or 
$C_s = q^{-2}$, $C_t = q$. Therefore we have respectively
$$
S {}^{ae}{}_{df} g^{fg} S {}^{cb}{}_{eg} 
= q^{\pm 2} g^{ac} \delta^b_d.
$$
In both cases $\sigma$ fulfills the braid Equation~(\ref{braid1}).  If
we compare the above equation with (\ref{2.2.24}) we see that the
metric~(\ref{fine`}) is not compatible with the linear connection. We
shall show however below in Equation~(\ref{line-element}) that with a
conformal factor it is so; the (flat) linear connection is equal to the
Levi-Civita connection of a (flat) metric conformally equivalent to the
one which we have found.  The symmetry condition (\ref{symm}) for the
metric follows from (\ref{utile1}) and from (\ref{decompo}).  Since the
matrix $g^{ab}$ is not symmetric in the ordinary sense there can exist
no linear transformation $\theta^{\prime a} = \Lambda^a_b \theta^b$ such
that $g^{\prime ab} = \delta^{ab}$.

{From} Equations~(\ref{ok}), (\ref{basic1}), and (\ref{gigio2}) it is
possible to determine the action of $\sigma$ and $g$ on the
basis $\xi^i$. One finds that
\be
\sigma (\xi^i\otimes \xi^j) =
S^{ij}{}_{hk}\,\xi^h\otimes \xi^k                            \label{gigio0}
\ee
and 
\be
g(\xi^i\otimes \xi^j) = 
g^{ij}\, \alpha^2 q^{-1}\,r^2\Lambda^2.                       \label{final}
\ee
Equation~(\ref{gigio0}) his is the same formula as (\ref{gigio2});
Equation~(\ref{final}) is to be compared with the
Equation~(\ref{fine`}).  From these results it is manifest that
$\sigma$ and $g$ are $SO_q(3)$-covariant; under the $SO_q(3)$ coaction
$\sigma(\xi^i\otimes \xi^j)$ and $g(\xi^i\otimes \xi^j)$ transform as
$\xi^i\otimes \xi^j$.  Relations (\ref{gigio0}), (\ref{final}) of course
could have been obtained also by a direct calculation in the $\xi^i$
basis.

For either choice of $\sigma$ the covariant derivative
\be
D_{(0)} \xi = 
-\theta \otimes \xi + \sigma (\xi \otimes \theta),            \label{conn}
\ee
defined in (\ref{2.2.14}) is manifestly $SO_q(3)$-invariant, because 
$\theta$ is. The most general $D$ is given by
Equation~(\ref{2.2.15}). But because of Equations~(\ref{ththcr}) we have
$$
\pi\circ \chi(\theta^a)  = - \chi^a{}_{bc} \theta^b \theta^c =
- \chi^a{}_{bc} (\c{P}_s +\c{P}_a +\c{P}_t)^{bc}_{de} \theta^d \theta^e =
- \chi^a{}_{bc} \c{P}_a{}^{bc}_{de} \theta^d \theta^e
$$
and because of Equations~(\ref{whynot}) this extra term must vanish if 
we require the torsion to vanish. Finally, one can show that a term 
$\chi^a{}_{bc} (\c{P}_s +\c{P}_t)^{bc}_{de} \theta^d \otimes\theta^e$
is forbidden by $SO_q(3)$-covariance. The most general torsion-free and 
$SO_q(3)$-covariant linear connection is given then by
\be
D = D_{(0)}.                                         \label{g-conn}
\ee

The linear curvature tensor is given by Equation~(\ref{2.20}); the first
term is absent, because $\theta^2 = 0$.  The result that the curvature
map $\mbox{Curv}$ is left-linear but in general not right-linear is
particularly evident from this formula.  In fact it is easy to see that
$\mbox{Curv} = 0$.  This relies only on the fact that $S$ fulfills the
braid equation (\ref{braid1}):
\bea
\mbox{Curv}(\xi) \squeeze&& = \pi_{12}\sigma_0{}_{12}
\sigma_0{}_{23}\sigma_0{}_{12}
(\theta^a\otimes \theta^b\otimes \theta^c)\:\xi_a\lambda_b\lambda_c
\nonumber\\[4pt]
&& = \pi_{12}(S_{12}S_{23}S_{12})^{abc}{}_{def}
(\theta^d\otimes \theta^e\otimes \theta^f)\:\xi_a\lambda_b\lambda_c
\nonumber\\[4pt]
&& = (S_{12}S_{23}S_{12})^{abc}{}_{def}
(\theta^d\theta^e\otimes \theta^f)\:\xi_a\lambda_b\lambda_c
\nonumber\\[4pt]
&& = -(S_{12}S_{23}\c{P}_a{}_{12})^{abc}{}_{def}
(\theta^d\theta^e\otimes \theta^f)\:\xi_a\lambda_b\lambda_c
\nonumber\\[4pt]
&& =  -(\c{P}_a{}_{23}S_{12}S_{23})^{abc}{}_{def}
(\theta^d\theta^e\otimes \theta^f)\:\xi_a\lambda_b\lambda_c
\nonumber\\[4pt]
&& = 0. \nonumber
\eea
The last three equalities follow from respectively
Equations~(\ref{2.2.6}), (\ref{braid2}) and (\ref{lambdacr1}). 

It is interesting to compute the result of the two flat connections on
the differentials. On one hand, if we choose $S = (q\hat R)^{-1}$ then
using Formulae~(\ref{xxicr}) and (\ref{utile1}) it is straightforward to
show that
\be
D\xi^i=0.                                                  \label{cicca}
\ee
This means that the `coordinates' are adapted to the zero-curvature
condition. On the other hand, it is straightforward to show that if 
$S = q\hat R$ then
\bea 
D\xi^i \squeeze&&=
(q^2-1)\theta \otimes \xi^i -(1+q^{-3})r^{-2}x^i [q^2\xi^l\otimes \xi^m
g_{lm}-q^4 g_{hj}\hat R^{ji}_{lm}x^h \xi^l\otimes \xi^m] \nonumber\\[4pt]
&&= (q^2-1)[\theta \otimes \xi^i +q^{-2}\xi^i\otimes \theta]
-q^2(1+q^{-1})r^{-2}x^i\xi^l\otimes \xi^mg_{lm}.            \label{explicit} 
\eea
This is different from zero. In this case the `coordinates' are not well
adapted. 

We can repeat the same construction for the barred calculus
$\bar\Omega^1(\c{A})$.  We simply list the results. There are
essentially a unique metric $g$ and two generalized permutations
$\sigma$ which fulfill the metric compatibility weakened by a conformal
factor. They are defined by
\bea
&&g(\bar\theta^a\otimes\bar\theta^b) = g^{ab}        \label{bbfine}\\[2pt]
&&\sigma (\bar\theta^a\otimes \bar\theta^b) =
\bar S^{ab}{}_{cd}\,\bar\theta^c\otimes \bar\theta^d,      \label{bbgigio2}
\eea
with either 
\be
\bar S = q\hat R, \qquad\mbox{or}\qquad \bar S = (q\hat R)^{-1}.
                                                      \label{bbdouble}
\ee
Correspondingly
$$
\bar S {}^{ae}{}_{df} g^{fg} \bar S {}^{cb}{}_{eg} 
= q^{\pm 2} g^{ac} \delta^b_d.
$$
In the $\bar\xi^i$ basis,
\bea
&&\sigma (\bar \xi^i\otimes \bar \xi^j) =
\bar S^{ij}{}_{hk}\,\bar\xi^h\otimes \bar\xi^k,  \label{bbgigio0} \\[2pt]
&&g(\bar \xi^i\otimes \bar\xi^j) = 
g^{ij}\, \alpha^2 q\, r^2\Lambda^{-2}.
                                                         \label{bbfinal}
\eea
The most general torsion-free $SO_q(3)$-covariant and (up to a conformal
factor) metric compatible covariant derivatives are
\be
\bar D \bar \xi = -\bar \theta \otimes \bar \xi + 
\sigma (\bar \xi \otimes \bar \theta)                          \label{bconn}
\ee
with either choice of $\sigma$.  They are manifestly
$SO_q(3)$-invariant.  The associated linear curvatures
$\overline{\mbox{Curv}}$ vanish.

If one extends the involution to the second tensor power of the 1-form
algebra using (\ref{TPI}), then it is easy to verify that neither $D$
nor $\bar D$ are real.  If $\bar S = S^{-1}$, they fulfill however
instead the relation
\be
(D \xi)^* = \bar D \xi^*.                                     \label{Dstar}
\ee

\sect{The commutative limit}                                  \label{climit}

The set $(x^i, e_a)$ generates the algebra $\c{D}_h$ of observables of
phase space; the $x^i$ generate the subalgebra of position observables,
whereas the $e_a$ the subalgebra of momenta observables. One has to look
for a complete set of three (the dimension of the classical underlying
manifold) independent commuting observables within the whole $\c{D}_h$,
and not just within the position space subalgebra, since the latter is
not abelian. Of course this is exactly what makes noncommutative
geometry interesting, in that it makes completely localized states
impossible.  Here we shall restrict our consideration especially to
position space observables.  It should be the theory which indicates
what is the correct identification of the operators $x^i$ in commutative
limit.  By theory we mean the algebra, that is, the theory proper, as
well as the representation or the state.  Below we find two possible and
rather different identifications of the $x^i$.  A physically
satisfactory understanding should involve a Hilbert space representation
theory of the algebra such that the spectrum of the position observables
becomes `dense' in the limit manifold.  Since the present analysis is
meant to be preliminary and heuristic in character, we shall 
restrict our attention to the algebra.

In the commutative limit $\Lambda$ commutes both with coordinates and
derivations, and therefore must be equal to a constant. Since we have
already normalized $\Lambda$ so as to be unitary, the constant must be a
pure phase. One can always absorb the latter in the same normalization,
so that
\be
\lim_{q\to 1} \Lambda =1.                                   \label{L-limit}
\ee 

In the first identification we pick a fixed real $\alpha$ in 
Equation~(\ref{stehbein3}) and suppose that the generators of the
algebra tend to their naive natural limit as (complex) coordinates on 
a real manifold and that the frame tends to the corresponding limit of
a moving frame on this manifold.
{From} (\ref{final}) we find that in the commutative limit the metric is
given by the line element \be ds^2 = (\alpha r)^{-2} \delta_{i,-j} dx^i
dx^j = (\alpha r)^{-2} (dx^2 + dy^2 + dz^2).  \label{line-element} \ee
The second expression is in terms of the real coordinates
(\ref{realmat}).  If one uses spherical polar coordinates then one sees
immediately that the Riemannian space is $S^2 \times \b{R}$ with 
$\log (\alpha r)$ the preferred coordinate along the line.  The radius
of the sphere is equal to $\alpha^{-1}$.  An interesting feature of this
identification is that neither the covariant derivative (\ref{g-conn})
nor (\ref{bconn}) is compatible with this non-flat metric.  This was of
course to be expected since none of the $\sigma$ we have used satisfies
the metric compatibility condition (\ref{2.2.24}). The problem we are
considering here lies in fact a little outside the range of the general
theory of Section~2 because of the element $\Lambda$ which is not in the
center but which has nevertheless a vanishing differential.
Alternatively if $\alpha$ is proportional rather to $r^{-1}$ the
conformal factor in (\ref{line-element}) would disappear and the metric
would be the ordinary flat metric of $\b{R}^3$.  We observed in
Section~4 that this would follow from the commutation relation
(\ref{zero'}).
In the commutative limit the frame~(\ref{stehbein3})
becomes a moving frame in the sense of Cartan. The singularity at
$x^0 = 0$ is to be expected , since there can be no global frame on a
non-parallelizable manifold like a sphere; for positive and negative
values of $x^0$ one has two different local sections of the frame bundle
on the two charts corresponding respectively to the upper and lower
hemisphere.  According to (\ref{bar-theta}), the frame is however not
real.  To find a real frame we first try the same linear transformation
(\ref{realmat}) used for the coordinates
$$
\theta^1 = {1 \over \sqrt 2} (\theta^- + \theta^+), \qquad
\theta^2 = \theta^0, \qquad
\theta^3 = {i \over \sqrt 2} (\theta^- - \theta^+).
$$
A short calculation yields
\bea
&&\theta^1 = (\alpha x^0 r)^{-1} (rdx - xdr + i z dr), \nonumber\\[2pt]
&&\theta^2 = (\alpha x^0 r)^{-1} (rdr - i x dz + i z dx), \label{clf} \\[2pt]
&&\theta^3 = (\alpha x^0 r)^{-1} (rdz - i xdr - z dr).
                                                             \nonumber
\eea
Although in the commutative limit the differential is real, we see that
the frame is not. Equivalently, from (\ref{bar-theta}) we see that the
frames for the two differential calculi do not coincide in the
commutative limit, although the two differential calculi themselves do.
It is often found in specific calculations in General Relativity that it
is more convenient to use a complex frame to calculate real curvature
invariants. We can see however no property of the frame (\ref{clf})
which makes it particularly adapted to study the space $S^2 \times \b{R}$,
or the space $\b{R}^3$ in the case that $\alpha\propto r^{-1}$.
A more complete analysis should involve at this point the study of the
$*$-representa\-tions of the algebra of observables $\c{D}_h$, studied
e.g. in Ref.~\cite{fioijmpa}. We shall not enter it here, but simply
recall that the definition of the commutative limit of these
representations is rather delicate. 

We shall now propose a different identification of the $x^i$ in the
commutative limit.  The limit manifold is the desired $\b{R}^3$ and
for finite $h$ the $x^i$ will be a sort of general coordinates on
$\b{R}^3$. They are related to cartesian coordinates by a
transformation which becomes singular in the commutative limit.  The
singularity can be removed by performing a renormalization procedure
(with diverging renormalization constant $\alpha^{-1}$
\footnote{Or with
a function $\alpha^{-1}=\alpha'r$, with $\alpha'$ a diverging constant,
in the case (\ref{zero'})})
before taking 
the limit. We thus obtain cartesian coordinates $y^i$ on $\b{R}^3$.
We try to solve algebraically the problems arising from the first
identification. We reconsider the search of a real nondegenerate frame
in the commutative limit.  If the condition 
$\lim\limits_{h\rightarrow 0}x=0=\lim\limits_{h\rightarrow 0} z$, or
equivalently
\be
x^{\pm}_{(0)}:=\lim\limits_{h\rightarrow 0}x^{\pm}=0,       \label{problem2}
\ee
were fulfilled on all points of the classical manifold, the frame
(\ref{clf}) would be real in the limit; apparently the latter would
also be degenerate ($\theta^1 = 0 = \theta^3$), but the latter
consequence can be avoided by choosing a diverging normalization
factor $\alpha^{-1}$ in (\ref{stehbein3}) and related definitions.

Let us first show by a one-dimensional example how such a
renormalization may change the situation.  In the sequel we shall append
a suffix ${(0)}$ to any quantity to denote its commutative limit, for
example $r_{(0)}=\lim_{h\rightarrow 0} r$.  Let $y\in\b{R}$ and
consider the change of coordinates
\be
x=f(\alpha y)
\ee
with some $\alpha$ going to zero when $h$ does and $f$ invertible, which
we shall assume normalized in such a way that $f'(0)=1$.  Although
$x_{(0)} = f(0) = \mbox{constant}$, for example, $f(0)=0$ we have
$$
dx_{(0)}=0, \qquad \frac{\partial}{\partial x}|_{h=0}=\infty.
$$
Nevertheless one can extract the original nontrivial coordinate,
differential and derivative by taking a difference and performing a
rescaling before taking the limit:
\be 
y=\lim\limits_{h\rightarrow 0}\frac{x-x_{(0)}}{\alpha} \qquad
\theta^1\equiv dy=\lim\limits_{h\rightarrow
0}\frac{dx-dx_{(0)}}{\alpha}\qquad 
\frac{\partial}{\partial y}=\lim\limits_{h\rightarrow 0} 
\alpha \frac{\partial}{\partial x}.                        \label{rescale} 
\ee 

Note that the above convergences are not uniform in $y$, namely are
slower as $y$ gets larger.  This example suggests that we may solve the
degeneracy problem by taking the normalization factor $\alpha$ in
formula (\ref{stehbein3}) as a suitable infinitesimal in $h$.

Now we return to the the algebra which interests us here. We aim to show
that, upon assuming (\ref{problem2}), and by choosing an infinitesimal
$\alpha$ such that the commutative limit of $\alpha^{-1} x^{\pm}$ is
finite, the commutative limit of $\alpha^{-1}
x^-,\alpha^{-1}(x^0-1),\alpha^{-1}x^+$ is in fact a set of (complex)
cartesian coordinates $y^-,y^0,y^+$ of $\b{R}^3$.  In addition, we shall
also sketch some basic ideas for a correspondence principle between the
deformed and the undeformed theory for finite $h$.  
As a first step, we show that one
can construct objects $y^i$ and $\frac{\partial}{\partial y^i}$, having
the commutation relations of a set of classical coordinates and their
derivatives, as `functions' of $x^i, e_a, h$ and a free parameter
$\alpha$. The `functions' are in fact power series in $h$. The zero
degree term for $y^i$ is essentially fixed by (\ref{problem2}).  We
stress in advance that the manipulations reported below are heuristic
and formal; these manipulations will have to be justified in the proper
representation theory framework.

The algebra of observables $\c{D}_h$ is a deformed Heisenberg algebra of
dimension 3 ($h$ is the deformation parameter). In fact, clearly the 
commutation relations between the $x^i$ are such that
\[
[x^i,x^j]=O(h).                            
\]
Moreover, the change of generators (\ref{covder})
is independent of $\alpha$ (the dependences in $\theta^a_i$
and $\lambda_a$ cancel), invertible, and gives
\bea
&&[\partial_i,x^j]=\delta^i_j +O(h), \nonumber\\
&&[\partial_i,\partial_j]=O(h).\nonumber
\eea
It is known (see e.g. ~\cite{pillin} and
references therein) that for any deformation $\c{D}_h$ of a Heisenberg
algebra $\c{D}$ there exists an isomorphism
$\c{D}_h\leftrightarrow\c{D}[[h]]$ of algebras over $\b{C}[[h]]$.
Applied to the present case, this implies that there exists a formal
realization of the canonical generators of $\c{D}$ in terms of
$x^i,\partial_j$, namely a set of formal power series in $h$
\bea
&& y^i=f^i(x,\partial,h)= x^i+ O(h)\nonumber\\[2pt]
&& \frac{\partial}{\partial y^i},
=g_i(x,\partial,h)=\partial_i+ O(h),\nonumber
\eea
with coefficients equal to polynomials in $x^i,\partial_j$ and
reducing to the identity for $h=0$, such that the objects
$y^i, \frac{\partial}{\partial y^i}$ fulfill the canonical
commutation relations
\be
[y^i,y^j]=0, \qquad 
[\frac{\partial}{\partial y^i},\frac{\partial}{\partial y^j}]=0
\qquad[\frac{\partial}{\partial y^i},y^j]=\delta^i_j.         \label{ccr}
\ee
When $h=0$ the above transformation is the identity, so also for 
$h\neq 0$ one can formally invert it to obtain a realization of the
`deformed generators' $x^i,\partial_i$ in the form of formal power
series in $h$ with coefficients equal to polynomials in
$y^i,\frac{\partial}{\partial y^i}$.  (As a consequence, also $e_a$ can
be expressed in the same form; one can show the same also for
$\Lambda$).  The above transformation is not unique. 
It is defined up an infinitesimal inner automorphism of the algebra
\be
x^i\rightarrow u\, x^i\, u^{-1} \qquad
\partial_i\rightarrow u\, 
\partial_i\, u^{-1},                 \label{innerauto}
\ee 
where $u=1+O(h)\in\c{D}_h[[h]]$. Since the commutation
relations (\ref{ccr}) are preserved upon
rescaling and translations of the $y^i$ by some constants 
$\alpha,c^i\in\b{C}$, $y^i\rightarrow c^i+\alpha y^i$,
if we relax the condition that the transformation
be the identity at zero order in $h$
a larger family of realizations will be given by
\bea
&& y^i= 
\alpha^{-1}[f^i(x,\partial,h)- c^i]=
\alpha^{-1}[x^i-c^i+ O(h)]\label{puppi1} \\
&& \frac{\partial}{\partial y^i}=\alpha\, g_i(x,\partial,h)=
\alpha[\partial_i+ O(h)].
\label{puppi2}
\eea
Next we shall fix the parameters $\alpha, c^i$ through some basic
requirements on the commutative limit. As for
the remaining indeterminacy (\ref{innerauto}), 
it should be fixed by geometrical requirements for
finite (though small) $h$, e.g. the ones suggested at the end of 
the section.

In the first identification, namely (\ref{line-element}),
we have picked $c^i=0$ and $\alpha$ finite,
e.g. $\alpha=1$; consequently $y^i=\lim_{h\to 0} x^i$. 
In the second identification we choose $c^{\pm}=0$ and $\alpha$ a suitable
infinitesimal in $h$ so that the limit (\ref{problem2}) is fulfilled;
consequently $y^{\pm}$ will be recovered as the limit 
$\lim_{h\to0} (x^{\pm}/\alpha)$. In order that formulae 
(\ref{puppi1}-\ref{puppi2})
contain a unique expansion parameter $\alpha$, instead of two
$\alpha,h$, we choose the infinitesimal as 
$\alpha=O(h^{1/p})$, with some $p\in\b{N}$ to be
determined. We also choose $c^0\neq 0$ 
in order that $x^0$ does not 
become singular in the commutative
limit. If for simplicity we normalize $c^0$
to 1, we will find $x^0_{(0)}=1$. Together with
(\ref{problem2}) and (\ref{squarelenght}), this implies
$r_{(0)}=1$. 
Note that any polynomial (or power series) of the $x^i,r$ with 
finite coefficients has a constant (i.e. independent of $y^i$)
commutative limit,  because $\alpha\to 0$.
We shall formally identify the three $y^i$ thus obtained as a
set of (local) classical coordinates on the limit manifold, and
$\frac{\partial}{\partial y^i}$ as the derivatives with respect to
$y^i$. Thus the limit manifold will have dimension three.
In the commutative limit 
\be
r\rightarrow 1,\qquad  (\alpha r)^{-1} dx^i\rightarrow dy^i, 
\ee
and the line element (\ref{line-element}) becomes
\be
ds^2  \to \sum_idy^i  dy^{-i}.          \label{limit-line-element}
\ee
This yields a vanishing Riemannian curvature and is consistent with
the vanishing (for all $q$) linear curvature computed in
Sect. \ref{metrics}.  Assuming a trivial topology, the limit
Riemannian manifold will be thus Euclidean space, and $y^i$ will be
(complex) cartesian coordinates on $\b{R}^3$.

{From} $r_{(0)}=x_{(0)}=\Lambda_{(0)}=1$ one also easily finds
$$
\lim\limits_{\alpha\rightarrow 0}
\alpha^{-1}q\Lambda e^i_a=\delta^i_a, \qquad
\lim\limits_{\alpha\rightarrow 0}
\alpha\Lambda^{-1} \theta^a_i=\delta^a_i,
$$
whence
\bea
&&\lim\limits_{\alpha\rightarrow 0} e_a
=\lim\limits_{\alpha\rightarrow 0}q\Lambda e^i_a\partial_i=
\lim\limits_{\alpha\rightarrow 0}\alpha^{-1}q\Lambda e^i_a
\frac{\partial}{\partial y^i}=\frac{\partial}{\partial y^a},\\[4pt]
&&\lim\limits_{\alpha\rightarrow 0} \theta^a
=\lim\limits_{\alpha\rightarrow 0}\Lambda^{-1} \theta^a_i dx^i=
\lim\limits_{\alpha\rightarrow 0}\Lambda^{-1} \theta^a_i
\alpha\alpha^{-1}dx^i=d y^a.                                \label{dopo}
\eea
This is consistent with the metric, stating that the set $\{\theta^a\}$
(and, in the dual formulation, the set $\{e_a\}$) is orthonormal.  Let
us determine possible values for $\alpha$. Plugging (\ref{puppi1}) in
the commutation relations (\ref{xxcr3}), it is easy to see that at
lowest order in $\alpha$ the third relation can be satisfied only if
$p\ge 3$. We shall later exhibit a transformation
(\ref{puppi1}-\ref{puppi2}) based indeed on $\alpha=O(h^{1\over 3})$.

To understand the physical meaning of $\alpha$ one would have to
consider the $*$-represen\-tations of $\c{D}$.  We shall see elsewhere
that it is related to the Planck lenght.  Moreover in the commutative
limit $y^0$ will `inherit' from the $*$-representation studied
previouly~\cite{fioijmpa} a spectrum dense in the whole $\b{R}$, in
agreement with the fact that $y^0$ becomes a real cartesian coordinate.
Similarly, $y_{\bot}^2:=y^+y^-$ `inherits' a spectrum dense in the whole
$\b{R}^+$, in agreement with the fact that it becomes the square
distance from the $y^0$ axis.

So far we have left the realization (\ref{puppi1}-\ref{puppi2}) still
largely undetermined by the inner automorphism freedom
(\ref{innerauto}). We may restrict this freedom by imposing further
requirements on $x^i,e_a$, so that the physical meaning of the latter 
at the
representation theoretic level become more manifest; for instance, since
they commute and are real, we may require that at all orders in $\alpha$
$x^0,r$ depend only on $\vec{y}$, so that their meaning of position
observables be manifest.  This will be discussed in detail
elsewhere. Once determined, the realization (\ref{puppi1}-\ref{puppi2})
will be an essential ingredient of the {\it correspondence principle}
between the deformed and undeformed theory 
for finite $h$. As an example, a formal realization of the
algebra fulfilling all previous requirements and such that the classical
involution $(y^i)^*=y^{-i}$ realizes also the involution $*$ of $\c{A}$
is
\bea
&& x^0= e^{\alpha y^0-\frac{\alpha^2}2} \label{special1} \\
&& \Lambda= e^{-\alpha^2\frac{\partial}{\partial y^0}}\\
&& x^+=\sqrt{\frac{e^{\alpha^2 y^+y^-}-1}{q^2(q+1)}}
e^{\alpha y^0 +i\varphi} e^{-\alpha^2\frac{\partial}{\partial y^0}
+ 2\alpha\frac{\partial}{\partial (y^+ y^-)}} \\
&& x^- =(x^+)^*\sqrt{q}
\eea
with $e^{\alpha^3}=q$ and $e^{2i\varphi}=y^+/y^-$. 
Its derivation will be given elsewhere. 
As a consequence
\be
r^2=e^{-2\alpha^3+\alpha^2y^+y^-+2\alpha y^0} \label{special2}
\qquad
(\frac{x^0}{r})^2=e^{\alpha^3-\alpha^2y^+y^-}.
\ee
Thus, the surfaces $r=const$ are paraboloids with axis $y^0$, the
surfaces $x^0=const$ are planes perpendicular to $y^0$, the surfaces
$x^0/r=const$ are cylinders with axis $y^0$, and the lines $x^0=const$,
$r=const$ are circles perpendicular to and with center on the axis
$y^0$.  
The exponential relation between $x^0$ and $y^0$ is analogous to the one
found \cite{CerHinMadWes98} for a 1-dimensional $q$-deformed model. From
the spectrum of $x^0$ \cite{fioijmpa} $y^0$ will inherit equidistant 
eigenvalues of both signs, suggesting a `uniform' structure of space. 
This will solve the problems \cite{fioijmpa} arising from the physical
identification of $x^i$ as cartesian coordinates.  
Summing up, the resulting geometry is now flat
$\b{R}^3$, instead of the $S^2\times \b{R}$ found in the first
identification (the one with finite $\alpha$).  $x^0, r$ 
may be adopted together with $\varphi$
as global curvilinear coordinates, but 
only for finite $\alpha\approx h^{1/p}$,
when they are to be identified with the coordinates (\ref{special1}) and
(\ref{special2}); for $\alpha\to 0$ they become singular, and coordinates
$y^i$ may be used at their place.

\sect{The involution and the real calculus}               \label{involu}

Formula (\ref{barstern}) gives an involution of 
$\Omega^1(\c{A})\oplus\bar\Omega^1(\c{A})$. Unfortunately, the latter
has rank 6 as a $\c{A}$-bimodule, instead of 3, and is not generated
from $\c{A}$ through the action of a real differential. 

The problem of constructing a differential calculus of rank 3 closed
under involution has
been considered by Ogievetsky and Zumino~\cite{olezu}, who expressed
$\bar\xi^i$ as functions not only of $x^i$, $\xi^i$, but also of
suitably $q$-deformed derivations $\partial_i$. A much more economical 
construction is the following. Note that one can define an algebra
isomorphism $\varphi:\Omega^*(\c{A})\rightarrow \bar\Omega^*(\c{A})$
acting as the identity on $\c{A}$ and on the 1-forms through
\be
\varphi(\theta^a)=\bar\theta^a,
\ee
since the commutation relations among the $\theta^a$'s 
and the $\bar\theta^a$'s are the same. Hence 
$\star= *\circ \varphi$ is an involution of $\Omega^*(\c{A})$
acting as
\be
(\theta^a)^{\star} = \theta^b g_{ba}.
\ee
This implies that $(\xi^i)^{\star}$ can be expressed as combinations
of $\xi^j$ with coefficients in (the extended) 
$\c{A}$. One can easily check that $\star$ has the
correct classical limit  provided the commutative
limit fulfills (\ref{problem2}).
Nevertheless, since $\varphi(d)\neq \bar d$, 
of course $d$ is not real under the involution.

It would be natural to define a real exterior derivative by
\[
d_r := (d + \bar d)  
\]
and a rank 3 $\c{A}$-bimodule 
$\Omega_r^1(\c{A})\subset\Omega^1(\c{A})\oplus\bar\Omega^1(\c{A})$ 
closed under $*$ from the generators
\[
\xi^i_r=\xi^i+\bar\xi^i.
\]
But this is impossible since
$\xi^i_r$'s do not close commutation relations with 
the $x^j$'s, as evident from (\ref{xxicr}), (\ref{xxicr'}).
To generate the exterior algebra through a real differential
we need to double either the number of the coordinates or 
the number of the 1-forms.
 
In general, consider an algebra $\c{A}$ with involution over which
there are two differential calculi $(\Omega^*(\c{A}), d)$ and
$(\bar\Omega^*(\c{A}), \bar d)$ neither of which is necessarily
real. Consider the product algebra $\tilde \c{A} = \c{A} \times \c{A}$
and over $\tilde \c{A}$ the differential calculus
\be
\tilde \Omega^*(\tilde \c{A}) =
\Omega^*(\c{A}) \times \bar\Omega^*(\c{A}).              \label{prod-cal} 
\ee 
It has a natural differential given by $\tilde d = (d, \bar d)$. The 
embedding
$$
\c{A} \hookrightarrow \tilde \c{A}
$$
given by $f \mapsto (f, f)$ is well defined and compatible with the
involution 
\be
(f,g)^* = (g^*,f^*)                                       \label{alg-invol}
\ee
on $\tilde \c{A}$.  Suppose there exists a frame $\theta^a$ for
$\Omega^*(\c{A})$ and a frame $\bar \theta^a$ for $\bar
\Omega^*(\c{A})$ and a relation of the form (\ref{bar-theta}) between
them which extends the involution (\ref{alg-invol}).  We can define
a real module $\Omega_r^1(\c{A})$ to be the
$\c{A}$-bimodule generated by the frame
$$
\theta_r^a = (\theta^a, \bar\theta^b g_{ba}).            \label{realframe}
$$
This does not necessarily contain however the image of 
$d_r = (d, \bar d)$. In fact in the case $\c{A} = \b{R}_q^3$ if we
require this to be the case then $\Omega_r^1(\c{A})$ is an 
$\tilde \c{A}$-bimodule and $\Omega_r^*(\c{A}) = \tilde\Omega^*(\c{A})$.

Alternatively one can consider a larger `function' algebra
containing two copies of $\c{A}$, whose generators we denote by
$x^i,\bar x^i$, with cross commutation relations of the form
\be
x^i\bar x^j = \hat R^{\pm 1}{}{ij}_{hk} \bar x^h x^k.
\ee
We can define an involution by setting
\be
(x^i)^*=\bar x^j g_{ji}                                  \label{invol}
\ee
and introduce an $SO_q(3)$-covariant differential calculus 
$\Omega(\tilde \c{A})$ with a real differential $d_r$ by setting 
\be
d_r x^i = \xi^i,\qquad d_r\bar x^i = \bar\xi^i,
\ee
with the commutation relations
\be
\ba{ll}
x^i\xi^j = q\hat R^{ij}_{hk}\xi^hx^k, 
&\bar x^i\bar \xi^j = q^{-1}\hat R^{-1}{}^{ij}_{hk}\bar\xi^h\bar x^k,\\[6pt]
x^i\bar\xi^j = q^{-1}\hat R^{\pm 1}{}^{ij}_{hk}\bar \xi^hx^k,
&\bar x^i\xi^j = q\hat R^{\mp 1}{}^{ij}_{hk} \xi^h\bar x^k.
\ea
\ee
(In the last two Equations we have picked a specific normalization).
These relations are compatible with (\ref{invol}) and $(d_r f)^*=d_r f^*$.
The commutation relations among the $\xi^i$ and
the $\bar\xi^i$ follow and are the same relations (\ref{xixicr}), 
(\ref{barxibarxicr}) found in Section \ref{calculi}
for the analogous elements in the conjugate calculi,
whereas the cross relations are 
\be
\xi^i\bar\xi^j = -q^{-1}\hat R^{\pm 1}{}^{ij}_{hk}\bar \xi^h\xi^k.
                                                        \label{barxixicr}
\ee
Thus one ends up with a real calculus with 6 independent coordinates
and 6 independent 1-forms.
But one can check that in this scheme there is no frame.

Finally, if we just let both the conjugate differentials $d,\bar d$
act on one copy of $\c{A}$ we generate a $\c{A}$-bimodule
of rank 6 closed under the involution $*$. An
enlarged exterior algebra is generated by $\xi^i,\bar\xi^i$ and
the generators of $\c{A}$; we need to add to the commutation
relations (\ref{xxicr}), (\ref{xixicr}), (\ref{xxicr'}), (\ref{barxibarxicr})
compatible commutation relations between $\xi^i$ and $\bar\xi^j$.
It is easy to show that the latter can only be of the form
$\xi^i\bar\xi^j =\gamma \hat R^ {-1}{}^{ij}_{hk}\bar\xi^h\xi^k$. 
This is compatible also with this algebra being a $\Omega^*(\c{A})$
and a $\bar\Omega^*(\c{A})$ bimodule.
If we extend $d,\bar d$ by requiring that 
\be
d\bar\xi^i = 0, \qquad\qquad \bar d\xi^i=0,                \label{barnilpo'}
\ee
$\gamma$ will be fixed to be $-q^{-1}$, and we find again
(\ref{barxixicr}).
It is immediate to show that $\theta$, $\bar\theta$ are the
Dirac operators of $d,\bar d$ also in the enlarged algebra.
A direct and lengthy calculation\footnote{We have performed it using
the program for symbolical computations REDUCE.}  shows that the
$\theta^a_i$, $\bar\theta^a_j$ introduced in Sect. 
\ref{calculi} satisfy the commutation relations
\be
\hat R^{ab}_{cd}\, \bar\theta^d_j\theta^c_i= 
\theta^b_l\bar\theta^a_k\, \hat R^{kl}_{ij}.             \label{basic'}
\ee
As a consequence of (\ref{barxixicr}), (\ref{basic'}), 
in the enlarged exterior algebra we find
\be
\theta^a\bar\theta^b =-q\hat R^{-1}{}^{ab}_{cd}\, \bar\theta^c\theta^d .
                                                        \label{th-bthcr}
\ee
We can also readily extend the generalized permutation $\sigma$
and the metric $g$ of Section \ref{metrics}
to $\Omega^1(\c{A})\oplus\bar\Omega^1(\c{A})$. 
For simplicity we consider an Ansatz where $\sigma$
does not mix the $\theta^a$'s and
the $\bar\theta^a$'s:
\bea 
\sigma (\theta^a\otimes \bar\theta^b) & = &
 V^{ab}{}_{cd}\,\bar\theta^c\otimes \theta^d, \label{bgigio2}\\
\sigma (\bar\theta^a\otimes \theta^b) & = &
\bar V^{ab}{}_{cd}\,\theta^c\otimes \bar\theta^d. \nonumber
\eea
Imposing conditions (\ref{2.2.6}), (\ref{th-bthcr}) we find
the numerical matrices $V,\bar V$:
\be
V=q\hat R^{-1}\qquad \qquad\bar V=q^{-1}\hat R.
\ee
Together with (\ref{gigio2}), (\ref{bbgigio2}),  these relations
completely define $\sigma$. If we set
\bea
&& g(\theta^a\otimes\bar\theta^b) = g^{a b}, \\
&& g(\bar\theta^a\otimes\theta^b) = \bar g^{a b}, 
\eea
the metric compatibility condition (\ref{met-comp}) 
on  $\Omega^1(\c{A})\otimes\bar\Omega^1(\c{A})$ and
$\bar\Omega^1(\c{A})\otimes\Omega^1(\c{A})$ reads
\be
\ba{l}
\bar V {}^{ae}{}_{df} g^{f g} S {}^{cb}{}_{eg} 
= g^{a c} \delta^b_d, \\
\bar V {}^{ae}{}_{df} \bar g^{f g} S {}^{cb}{}_{eg} 
= \bar g^{ac} \delta^b_d.
\ea
\label{brrr}
\ee
As well as its barred counterpart,
it can be satisfied only if the numerical matrices 
$g^{ab}$, $\bar g^{ab}$ are proportional to the 
matrix  (\ref{metric}) and we restrict
(\ref{double}) to the choice
\be
S=q\hat R, \qquad\qquad  \bar S= (q\hat R)^{-1}.         \label{unique}
\ee
Note that in (\ref{brrr}) no conformal factor appears.
The connections defined by (\ref{conn}), 
(\ref{bconn}) are now automatically extended to 
the exterior algebra of rank 6, and it is easy to check that also
on this larger domain their linear curvatures are zero. Moreover, they
are still related by (\ref{Dstar}). We conclude by noting that
\[
ds^2:=\theta^a\otimes\bar\theta^b g_{ab} + \mbox{h. c.}
\]
is real and annihilated by the action of both $D,\bar D$
(in other words it is ``invariant under parallel transport'').
Therefore it 
might be considered as a candidate for the square displacement.

\sect{Conclusions}

The fundamental open problem of the noncommutative theory of gravity
concerns of course the relation it might have to a future quantum theory
of gravity either directly or via the theory of strings and membranes.
But there are more immediate technical problems which have not received
a satisfactory answer. The most important ones concern the definition of
the curvature. It is not certain that the ordinary definition of
curvature taken directly from differential geometry is the quantity
which is most useful in the noncommutative theory. The main interest of
curvature in the case of a smooth manifold definition of space-time is
the fact that it is local.  We have defined Riemann curvature as a map
Curv which takes $\Omega^1(\c{C}(V))$ into
$\Omega^2(\c{C}(V)) \otimes_{\c{C}(V)} \Omega^1(\c{C}(V))$.
If $\xi \in \Omega^1(\c{C}(V))$ then $\mbox{Curv}(\xi)$ at a given
point depends only on the value of $\xi$ at that point. This can be
expressed as a bilinearity condition; the above map is a 
$\c{C}(V)$-bimodule map.  If $f \in \c{C}(V)$ then
\be
f \mbox{Curv}(\xi) = \mbox{Curv}(f\xi), \qquad  
\mbox{Curv}(\xi f) = \mbox{Curv}(\xi) f.
\ee
One would like to insure also that one is dealing with the
noncommutative version of a real manifold. This reality condition
exchanges left and right linearity and adds weight to the argument that
a bilinear curvature map is necessary.  In the noncommutative case
bilinearity is therefore the natural and only possible expression of
locality. It has not yet been possible to enforce it in a satisfactory
manner~\cite{DubMadMasMou96}.

We have argued here in favour of one possible extension to
noncommutative geometry of the definition of a linear connection which
relies essentially on the above expression of reality and locality.
Alternative extensions have however been given. In fact one definition
has been proposed~\cite{CunQui95} which is indeed local in the sense we
have defined locality but which is valid only in the noncommutative case
and becomes singular in the commutative limit.  The main difference with
alternative definitions and the one we use lies in the fact that the
locality condition is relaxed and the Leibniz rule is applied only from
one side; the module of 1-forms is considered only as a left (or right)
module. This means in fact that gravity is considered as a Yang-Mills
field with the general linear group or some $q$-deformation of it, as
structure group. One of the authors has given elsewhere~\cite{Mad97} a
partial list of these alternative proposals. We mention here only three
which are especially relevant to the example of quantum Euclidean space.
Ref. \cite{HeckSchmue95} is the closest to ours. 
A quantum group ($SL_q(N),O_q(N)$ or $Sp_q(N)$) function algebra, 
rather than a quantum space function algebra, is
chosen as algebra $\c{A}$.  General definitions
of a metric, covariant derivative, torsion, curvature, metric-compatible
covariant derivative on $\c{A}$ are given; in particular, the quantum 
group invariant/bicovariant ones are determined. The definitions
differ from ours essentially in that 
$\c{A}$-bilinearity is weakened to a left- (or right-) linearity.
The existence of frames is not investigated. In Ref.~\cite{Cas95} a
$q$-Poincar\'e group function algebra is chosen as an algebra $\c{A}$; the
definition of differential calculus, linear connection, curvature, mimic
the group-geometric ones on the ordinary Poincar\'e group manifold. 
The underlying Hopf algebra is triangular, what essentially means 
$\hat R^2=1$ (in this respect the situation is simpler
than the one considered in the present work).  The existence of metric,
torsion and frames is not investigated. The $q$-Minkowski space,
its differential calculus, linear connection, curvature, are meant to be
derived by projecting out the $q$-Lorentz group.
Other 
authors~\cite{Dur98, maj97} have abandoned the definition of a
connection as a covariant derivative and proposed a $q$-deformed version
of a principle bundle to define a noncommutative `Riemannian geometry',
using for example the universal calculus.  It is interesting to note
that according to the definition we use here the unique linear
connection associated with any universal calculus is necessarily
trivial~\cite{DubMadMasMou95}; the covariant derivative coincides with
the ordinary exterior derivative.

Finally, we would like to compare the results found here for the quantum
Euclidean space with the results found in~\cite{DimMad96} for the Manin
quantum plane. We start by noting that all commutation relations in both
cases are homogeneous separately in $x^i$ and $\xi^i$.  A Dirac operator
$\theta$ (\ref{fund}) for a quantum group covariant differential
calculus must be of degree 1 in $\xi^i$, degree -1 in $x^j$ and quantum
group invariant. There can be no such object for the Manin plane
because, to mimic the construction (\ref{deftheta}), we would have to
use the isotropic tensor $\varepsilon_{ij}$ of $SL_q(2)$, the
$q$-deformed $\varepsilon$-tensor in the place of the isotropic tensor
$g_{ij}$ of $SO_q(3)$, but the analog $x^ix^j\varepsilon_{ij}$ of $r^2$
vanishes, by the $q$-commutation relations of the Manin plane. The
absence of such a $\theta$ prevents the construction of a torsion-free
linear connection by means of formula (\ref{2.2.14}) which, by means of
Equation~(\ref{2.2.19}), could eventually yield a bilinear zero
curvature.  In Reference~\cite{DimMad96} a Stehbein was constructed
using the inverses of both Manin plane coordinates and no dilatator,
whereas we have needed here the inverse of one coordinate $x^0$ and the
introduction of the dilatator. Using the Stehbein, also a metric has
been found. However, it is not a quantum group covariant metric, as
found here.  One could define a $SL_q(2)$-covariant metric $g_0$ on the
Manin plane by introducing a dilatator $\Lambda$, in the same way as in
Equation~(\ref{lambda}), (\ref{zero}) and by setting $g_0(\xi^i\otimes
\xi^j)=\Lambda^{-3}\varepsilon^{ij}$ (See Equation (\ref{final}).
This would be a symplectic form in the limit $q=1$. 
An analog of equation (\ref{brrr}), but in the $\xi^i$ basis, would hold.

\section*{Acknowledgment} 
One of the authors (JM) would like to thank the J. Wess for his
hospitality and the Max-Planck-Institut f\"ur Physik in M\"unchen for
financial support.


\begin{thebibliography}{99}
\parskip 5pt plus2pt minus2pt
\tolerance=1000

\bibitem{cawa} 
U.\ Carow-Watamura, M.\ Schlieker and S.\ Watamura, ``$SO_q(N)$ covariant
Differential Calculus on Quantum Space and Quantum 
Deformation of Schroedinger Equation'',
Z. Phys. C Part. Fields {\bf 49} (1991) 439.

\bibitem{Cas95}
L.\ Castellani, ``Differential Calculus on $ISO_q(N)$, Quantum Poincar\'e
algebra and $q$-Gravity, Commun. Math. Phys. {\bf 171} (1995) 383.

\bibitem{CelGiaSorTar91} 
E. Celeghini, R. Giachetti, E, Sorace, M. Tarlini, 
``The Three-dimensional Euclidean Quantum Group $E(3)-q$ and its
$R$-matrix'', J. Math. Phys. {\bf 32} (1991) 1159.

\bibitem{CerHinMadWes98} 
B.\ L.\ Cerchiai, R.\ Hinterding, J.\ Madore, J.\ Wess
``The Geometry of a $q$-deformed Phase Space'', 
Munich Preprint LMU 98/08, math.QA/9807123.

\bibitem{Con94}
A.\ Connes, ``Noncommutative Geometry'', 
Academic Press, 1994.

\bibitem{CunQui95}
A.\ Cuntz, D.\ Quillen, ``Algebra extensions and nonsingularity'',
J. Amer. Math. Soc. {\bf 8} (1995) 251.

\bibitem{DimMad96}
A.\ Dimakis, J.\ Madore, ``Differential Calculi and Linear Connections'', 
J. Math. Phys. {\bf 37} (1996) 4647.

\bibitem{DubMadMasMou95}
M.\ Dubois-Violette, J.\ Madore, T.\ Masson, J.\ Mourad, 
``Linear Connections on the Quantum Plane'', 
Lett. Math. Phys. {\bf 35} (1995) 351.

\bibitem{DubMadMasMou96}
M.\ Dubois-Violette, J.\ Madore, T.\ Masson, J.\ Mourad, 
``On Curvature in Noncommutative Geometry'',
J. Math. Phys. {\bf 37} (1996) 4089.

\bibitem{Dur98}
M.\ Durdevi\'c, ``Differential Structures on Quantum Principal
Bundles'', Rep. on Math. Phys. {\bf 41} (1998) 91.

\bibitem{frt}
L.\ D.\ Faddeev, N.\ Y.\ Reshetikhin and L.\ A.\ Takhtajan,
``Quantization of Lie Groups and Lie Algebras'', 
Algebra i Analysis, {\bf 1} (1989),
178; translation: Leningrad Math.\ J.\ {\bf 1} (1990), 193.

\bibitem{fiothesis} 
G.\ Fiore, ``q-Euclidean Covariant Quantum mechanics on 
$\b{R}^N_q$: Isotropic Harmonic Oscillator and Free Particle'',
Ph. D. Thesis, SISSA-ISAS (Trieste), 1994.

\bibitem{fiodet} 
G.\ Fiore, ``Quantum Groups  $SO_q(N)$, $Sp_q(n)$ 
have $q$-Determinant, too'', J. Phys. A: Math. Gen. {\bf 27} (1994), 1-8.

\bibitem{fioijmpa} G.\ Fiore, ``The $q$-Euclidean Algebra $U_q(e^N)$ and
the Corresponding $q$-Euclidean Lattice'', Int. J. Mod. Phys. {\bf A11}
(1996), 863-886.

\bibitem{FioMad98}
G.\ Fiore, J.\ Madore, ``Leibniz Rules and Reality Conditions'',
Naples Preprint 98-13, math/9806071.

\bibitem{HeckSchmue95}
I.\ Heckenberger, K.\ Schmuedgen,
``Levi-Civita Connections on the Quantum Groups $SL_q(N), O_q(N)$ and
$Sp_q(N)$'', Commun. Math. Phys. {\bf 185} (1997), 177.
     
     
\bibitem{Kos60}
J.L. Koszul, ``Lectures on Fibre Bundles and Differential Geometry'',
Tata Institute of Fundamental Research, 1960, Bombay.

\bibitem{Lan97}
G. Landi, ``An Introduction to Noncommutative Spaces and their Geometries'',
Springer Lecture Notes, Springer-Verlag, 1997.

\bibitem{Mad95}
J. Madore, ``An Introduction to Noncommutative Differential Geometry
and its Physical Applications'', Cambridge University Press (1995).

\bibitem{Mad97}
J. Madore, ``Gravity on Fuzzy Space-Time'',
Vienna Preprint ESI 478, gr-qc/9709002.

\bibitem{MadMou98}
J. Madore, J. Mourad, ``Quantum Space-Time and Classical Gravity'', 
J. Math. Phys. {\bf 39} (1998) 423.

\bibitem{maj2} S. Majid, ``Braided Momentum in the $q$-Poincar\'e group'',
J. Math. Phys. {\bf 34} (1993), 2045.
 
\bibitem{maj97}
S. Majid, ``Quantum and Braided Group Riemannian geometry'',
Preprint Damtp/97-73, q-alg/9709025.

%
\bibitem{oleg2}
O.\ Ogievetsky, 
``Differential operators on quantum spaces for $GL_q(n)$ and $SO_q(n)$''
Lett. Math. Phys. {\bf 24} (1992), 245.
%%

\bibitem{wessvari}
 O. Ogievetsky, W. B. Schmidke, J. Wess and B. Zumino, 
``$q$-deformed Poincar\'e algebra'', Commun. Math.
{\bf 150} (1992) 495-518.

\bibitem{olezu}
O. Ogievetsky and B. Zumino, ``Reality in the Differential Calculus on
$q$-Euclidean Spaces'', Lett. Math. Phys. {\bf 25} (1992) 121-130.
   
\bibitem{pillin} M.\ Pillin, 
``On the Deformability of Heisenberg Algebras'', 
Commun. Math. Phys. {\bf 180} (1996), 23.

\bibitem{munich}
M. Schlieker, W. Weich and R. Weixler, 
``Inhomogeneous Quantum Groups'', Z. Phys. {\bf C 53} (1992), 79-82.

\bibitem{Sny47}
H.S. Snyder, ``Quantized Space-Time'', 
Phys. Rev. {\bf 71} (1947) 38.

\bibitem{WesZum90}
J. Wess, B. Zumino,
``Covariant differential calculus on the quantum hyperplane'',
Nucl. Phys. (Proc. Suppl.) {\bf B18} (1990) 302.

\end{thebibliography}
\end{document}